\documentclass[a4paper,10pt]{article}
\usepackage{amsmath,amssymb,vmargin,bbm,mathrsfs,amsfonts,geometry,graphicx,url}
\usepackage{color}
\title{\textbf{Entrance and sojourn times for Markov chains. \\
Application to $(L,R)$-random walks}}
\author{Valentina \textsc{Cammarota}\footnote{Dipartimento di Matematica,
\textsc{Universit\`a degli Studi di Roma `Tor Vergata'},
Via della Ricerca Scientifica 1, 00133 Roma,
\textsc{Italy}. \textit{E-mail address}: \texttt{cammarot@mat.uniroma2.it},
\textit{Web page:} \texttt{https://sites.google.com/site/} \texttt{valentinacammarota}}
\ and Aim\'e \textsc{Lachal}\footnote{
P\^ole de Math\'ematiques/Institut Camille Jordan/CNRS UMR5208, B\^at. L. de Vinci,
\textsc{Institut National des Sciences Appliqu\'ees de Lyon},
20 av. A. Einstein, 69621 Villeurbanne Cedex, \textsc{France}.
\textit{E-mail address:} \texttt{aime.lachal@insa-lyon.fr},
\textit{Web page:} \texttt{http://maths.insa-lyon.fr/$\scriptstyle\sim$alachal}}
}
\date{}

\setmarginsrb{20mm}{10mm}{20mm}{10mm}{0mm}{10mm}{0mm}{10mm}

\newtheorem{teo}{Theorem}[section]
\newtheorem{prop}{Proposition}[section]

\newtheorem{lem}{Lemma}[section]
\newtheorem{cor}{Corollary}[section]
\newtheorem{rem}{Remark}[section]
\newtheorem{ex}{Example}[section]
\numberwithin{equation}{section}
\newcommand{\Dim}{\noindent\textsc{Proof\\}}
\newcommand{\lqn}[1]{\noalign{\noindent $\displaystyle{#1}$}}
\allowdisplaybreaks

\begin{document}
\maketitle

\begin{abstract}
\noindent In this paper, we provide a methodology for computing the probability
distribution of sojourn times for a wide class of Markov chains. Our methodology
consists in writing out linear systems and matrix equations for generating
functions involving relations with entrance times.
We apply the developed methodology to some classes of random walks with
bounded integer-valued jumps.
\end{abstract}

\begin{small}

\end{small}

\section{Introduction}

We first introduce some settings: we denote by $\mathbb{Z}$ the set of all
integers, by $\mathbb{Z}^+$ that of positive integers:
$\mathbb{Z}^+=\mathbb{Z}\cap (0,+\infty)=\{1,2,\dots\}$,
by $\mathbb{Z}^-$ that of negative integers:
$\mathbb{Z}^-=\mathbb{Z}\cap(-\infty,0)=\{\dots,-2,-1\}$, and
by $\mathbb{Z}^{\dag}$ that of non-negative integers:
$\mathbb{Z}^{\dag}=\mathbb{Z}^+\cup \{0\}=\{0,1,2,\dots\}$.
We have the partition $\mathbb{Z}=\mathbb{Z}^+\cup \{0\}\cup \mathbb{Z}^-$.
These settings will be related to space variables.
We introduce other settings for time variables: $\mathbb{N}$ is the set
of non negative integers: $\mathbb{N}=\{0,1,2,\dots\}$ and $\mathbb{N}^*$ is
that of positive integers: $\mathbb{N}^*=\{1,2,\dots\}$.

While nearest neighbour random walk on $\mathbb{Z}$ has been extensively
studied, random walk with several neighbours seems to be less considered.
In this paper, we consider this class of random walks, i.e., those
evolving on $\mathbb{Z}$ with jumps of size not greater than $R$
and not less than $-L$, for fixed positive integers $L,R$, (the so-called
$(L,R)$-random walks). More specifically, we are interested in the time spent
by the random walk in some subset of $\mathbb{Z}$, e.g., $\mathbb{Z}^{\dag}$
up to some fixed time.

\subsection{Nearest neighbour random walk}

Let us detail the most well-known example of random walk on $\mathbb{Z}$:
let $(S_m)_{m\in\mathbb{N}}$ be the classical symmetric Bernoulli random walk
defined on $\mathbb{Z}$. The probability distribution of the sojourn time of
the walk $(S_m)_{m\in\mathbb{N}}$ in $\mathbb{Z}^{\dag}$ up to a fixed time
$n\in\mathbb{N}^*$,
$$
T_n=\#\{m\in\{0,\dots,n\}:S_m\ge 0\}=\sum_{m=0}^n\mathbbm{1}_{\mathbb{Z}^{\dag}}(S_m),
$$
is well-known. A classical way for calculating it consists in using generating
functions; see, e.g., \cite[Chap. III, \S 4]{feller} or
\cite[Chap. 8, \S 11]{renyi} for the case of the nearest neighbour random walk
on $\mathbb{Z}$, and \cite[Chap. XIV, \S 8]{feller} where the case of random
walk with general integer-valued jumps is mentioned. The methodology consists in
writing out equations for the generating function of the family of numbers
$\mathbb{P}\{T_n=m\}$, $m,n\in\mathbb{N}$.
A representation for the probability distribution of $T_n$ can be
derived with the aid of Sparre Andersen's theorem (see~\cite{sparre1,sparre2}
and, e.g., \cite[Chap. IV, \S20]{spitzer}).
Moreover, rescaling the random walk and passing to the limit, we get the
most famous Paul Levy's arcsine law for Brownian motion.

Nevertheless, the result is not so simple.
By modifying slightly the counting process of the positive terms of the random
walk, as done by Chung \& Feller in \cite{chung}, an alternative sojourn time for
the walk $(S_m)_{m\in\mathbb{N}}$ in $\mathbb{Z}^{\dag}$ up to time~$n$ can be
defined:
\begin{equation}\label{modified-sojourn}
\tilde{T}_n=\sum_{m=1}^n\delta_m\text{ with }\delta_m=\begin{cases}
1 &\text{if $(S_m>0)$ or $(S_m=0$ and $S_{m-1}>0)$,}
\\
0 &\text{if $(S_m<0)$ or $(S_m=0$ and $S_{m-1}<0)$.}
\end{cases}
\end{equation}
In $\tilde{T}_n$, one counts each time $m$ such that $S_m>0$ and only those
times such that $S_m=0$ which correspond to a downstep: $S_{m-1}=1$.
This convention is described in~\cite{chung} and, as written
therein--\textit{``The elegance of the results to be
announced depends on this convention''} (\textit{sic})--, it produces a
remarkable result. Indeed, in this case, the probability distribution of this
sojourn time takes the following simple form: for even integers $m,n$ such that
$0\le m\le n$,
$$
\mathbb{P}\{\tilde{T}_n=m\}=\frac{1}{2^n}\binom{m}{m/2}\!\binom{(n-m)}{(n-m)/2}.
$$
In addition, also the conditioned random variable $(\tilde{T}_n|S_n=0)$, for
even~$n$, is very simple: it is the uniform distribution on $\{0,2,\dots,n-2,n\}.$
The random variables $(\tilde{T}_n|S_n>0)$ and  $(\tilde{T}_n|S_n<0)$
admit remarkable distributions as well; see, e.g., \cite{random-walk}.
Hence, it is interesting to work with the joint distribution of $(\tilde{T}_n,S_n)$.

We observe that this approach could be adapted to a wider range of
Markov chains.

\subsection{Main results}

In this paper, we consider a large class of Markov chains on
a finite or denumerable state space $\mathcal{E}$ and we introduce the
time~$T_n$ spent by such a chain in a fixed subset $E^{\dag}$ of $\mathcal{E}$
up to a fixed time~$n$. Inspired by the modified counting
process~(\ref{modified-sojourn}), we also introduce an alternative
sojourn time $\tilde{T}_n$.
Let us recall that sojourn times in a fixed state play a fundamental role in the framework of
general Markov chains and potential theory.

We develop a methodology for computing the probability distributions of
$T_n$ and $\tilde{T}_n$ via generating functions (Theorems~\ref{th-Kixy} and~\ref{th-Ktildeixy}).
The technique is the following: by applying the Markov property, we write out
linear systems for the generating functions of $T_n$ and $\tilde{T}_n$.
Though it seems difficult to solve these systems explicitly, nevertheless
they could be numerically solved.
We refer the reader to the famous book \cite[Chap. XVI]{feller},
to \cite{kemeny-snell}, \cite{kemeny-snell-knapp}, \cite{norris}, \cite{woess}
for an overview on Markov chains and to the recent book~\cite{flajolet}
for generating functions.

Next, we apply the general results we obtained to the case of $(L,R)$-random walk
(Theorems~\ref{th-Kixy-LR}, \ref{th-tilde-Kixy-LR} and \ref{th-K-Ktilde-LR}).
We exhibit explicit results in the particular cases where $L=R=1$, namely that
of nearest random walk with possible stagnation (Theorems~\ref{th-K-ordinary}
and~\ref{th-Ktilde-ordinary}), and where $L=R=2$, the case of the two-nearest
random walk (Theorem~\ref{th-gene-M=2}). In this latter, we illustrate the
matrix approach which can be completely carried out.

In the case of the usual Bernoulli random walk on $\mathbb{Z}$, the
determination of the generating function of the sojourn time in
$\mathbb{Z}^{\dag}$ goes through that of the first hitting time of~$0$.
In this particular example, level~$0$ acts as a boundary for
$\mathbb{Z}^{\dag}$ in the ordered set $\mathbb{Z}$.
This observation leads us naturally to define in our context of Markov chains
on $\mathcal{E}$ a kind of boundary, $E^{\rm o}$ say, for the set of
interest $E^{\dag}$ which is appropriate to our problem (see
Section~\ref{section-settings}). This is the reason why we restrict our study
to Markov chains satisfying Assumptions~$(A_1)$ and $(A_2)$ (see
Section~\ref{section-settings}).
Let us mention that, as in the case of Bernoulli walk, entrance times
play an important role in the analysis.
\\

\subsection{Background and motivation}

Our motivations come--among others--from biology.
Certain stochastic models of genomic sequences are based on random walks
with discrete bounded jumps. For instance, DNA and protein sequences,
which are made of nucleotides or amino acids, can be modelled by Markov chains
whose state space is a finite alphabet. Biological sequence analysis can be
performed thanks to a powerful indicator: the so-called local score which
is a functional of some random walk with discrete bounded jumps.
This indicator plays an important role for studying alignements of two
sequences, or for detecting particular functional properties of a sequence
by assessing the statistical significance of an observed local score;
see, e.g., \cite{mercier}, \cite{mercier2} and references therein.

Another example in biology concerns the study of micro-domains
on a plasmic membrane in a cellular medium. The plasmic membrane is a place
of interactions between the cell and its direct external environment.
A naive stochastic model consists in viewing several kinds of constituents
(the so-called ligands and receptors) as random walks evolving on the
membrane. The time that ligands and receptors bind during a fixed amount of
time plays an important role as a measurement of affinity/sentivity of ligands
for receptors. It corresponds to the sojourn time in a suitable set for a
certain random walk; see, e.g., \cite{cholesterol}.

Let us mention other fields of applications where entrance times, exit times,
sojourn times for various random walks are decisive variables: finance,
insurance, econometrics, reliability, management,
queues, telecommunications, epidemiology, population dynamics...

More theoretically, several authors considered $(L,R)$-random walks
(i.e., with integer-valued jumps lying between $-L$ and $R$ where $L$ and $R$
are two positive integers) in the context of random walks in random environment.
When fixing the environment, the quenched law they deal with is associated with
a $(L,R)$-random walk. Let us quote for instance a pioneer work \cite{key}; next
\cite{der1} and \cite{der2} where the jumps are~$+2$ and~$-1$ (viz.
$(1,2)$-random walk); \cite{bremont} and recent papers \cite{hong-wang1},
\cite{hong-wang2}, \cite{hong-zhang}, corresponding to the particular
cases $L=1$ or $R=1$, namely the jumps are~$1$ in one direction and greater
than $1$ in the other direction. Many other references can be found therein.

In all the aforementioned papers, one of the main motivation is the study of
recurrence/transience and the statement of a law of large numbers for random
walk in random environment. Important tools for tackling this study are passage
time, exit time from a bounded interval, ladder times, excursions... in the
quenched (i.e. fixed) environment. For instance, in \cite{hong-wang2}, the
authors study passage time and sojourn time above a level before reaching
another one. In \cite{hong-wang1}, the authors are able to compute generating
function of the exit time of $(2,2)$-random walk (see also \cite{bremont} for
the case of $(1,2)$-random walk). In particular, in \cite{hong-zhang}, we can read:
\textit{``These exit probabilities [in deterministic environment]
play an important role in the offspring distribution of the branching
structure, which can be expressed in terms of the environment'' (sic).}

\subsection{Plan of the paper}

The paper is organized as follows. In Section~\ref{section-settings} we
introduce the settings. In particular, inspired by~(\ref{modified-sojourn}),
we elaborate an alternative counting process.
In Section~\ref{section-gene-func} we consider several
generating functions and, in particular, in Theorems~\ref{th-Kixy} and~\ref{th-Ktildeixy},
we describe a method for computing the generating functions of $(T_n, X_n)$
and $(\tilde{T}_n,X_n)$ for a general Markov chain satisfying Assumptions~$(A_1)$ and $(A_2)$.
Since the proofs of Theorem~\ref{th-Kixy} and \ref{th-Ktildeixy} are quite
technical, we postpone them to Section~\ref{proof-section} as well
as those of auxiliary results. In Section~\ref{section-LR} we apply the
methodology developed for general Markov chains to the class of $(L,R)$-random
walks by adopting a matrix approach. Finally Section~\ref{subsection-ordinary-RW}
and Section~\ref{subsection-symmetric} are devoted to the more striking
examples of ordinary random walks and symmetric $(2,2)$-random walks.


\section{Settings}\label{section-settings}

Let $(X_m)_{m\in\mathbb{N}}$ be an homogeneous Markov chain on a state space $\mathcal{E}$
(which is assumed to be finite or denumerable) and let $E^{\dag}$, $E^{\rm o}$ be subsets
of $\mathcal{E}$ with $E^{\rm o}\subset E^{\dag}$. We set $E^+=E^{\dag}\backslash E^{\rm o}$,
$E^-=\mathcal{E}\backslash E^{\dag}$ and $E^{\pm}=\mathcal{E}\backslash E^{\rm o}=E^+\cup E^-$.
Actually, we partition the state space into $\mathcal{E}=E^+\cup E^{\rm o}\cup E^-$.
We will use the classical convention
that $\min(\varnothing)=+\infty$ and that $\sum_{\ell=p}^q=0$ if $p>q$.
Throughout the paper, the letters
$i,j,k$ will denote generic space variables in $\mathcal{E}$ and $\ell,m,n$
will denote generic time variables in $\mathbb{N}$. We also introduce the
classical conditional probabilities
$\mathbb{P}_i\{\cdots\}=\mathbb{P}\{\cdots | X_0=i\}$
and we set $p_{ij}=\mathbb{P}_i\{X_1=j\}$ for any states
$i,j\in\mathcal{E}$.

\subsection{Entrance times}

It will be convenient to introduce the first entrance times
$\tau^{\rm o},\tau^{\dag},\tau^+,\tau^-,\tau^{\pm}$
in $E^{\rm o},E^{\dag},E^+,E^-,E^{\pm}$ respectively:
\begin{align*}
\tau^{\rm o} =\, & \min\{m\in\mathbb{N^*}: X_m\in E^{\rm o}\},\\
\tau^{\dag}  =\, & \min\{m\in\mathbb{N^*}: X_m\in E^{\dag}\},\\
\tau^+       =\, & \min\{m\in\mathbb{N^*}: X_m\in E^+\},\\
\tau^-       =\, & \min\{m\in\mathbb{N^*}: X_m\in E^-\},\\
\tau^{\pm}   =\, & \min\{m\in\mathbb{N}^*: X_m\in E^{\pm}\}.
\end{align*}
We plainly have $\tau^{\dag}\le\tau^{\rm o}$ and $\tau^{\pm}=\tau^+\wedge\tau^-$.
Additionally, we will use the first entrance time in $E^{\pm}$ after time~$\tau^{\rm o}$:
$$
\tilde{\tau}^{\pm}=\min\{m\ge \tau^{\rm o}: X_m\in E^{\pm}\}.
$$
As a mnemonic, our settings write in the example of the classical Bernoulli random walk
on $\mathcal{E}=\mathbb{Z}$, with the choices $E^{\rm o}=\{0\}$,
$E^{\dag}=\mathbb{Z}^{\dag}=\{0,1,2,\dots\}$,
$E^+=\mathbb{Z}^+=\{1,2,\dots\}$, $E^-=\mathbb{Z}^-=\{\dots,-2,-1\}$ and
$E^{\pm}=\mathbb{Z}\backslash\{0\}$, as
\begin{align*}
\tau^{\rm o}=\, & \min\{m\in\mathbb{N}^*: X_m=0\},
\\
\tau^{\dag} =\, & \min\{m\in\mathbb{N}^*: X_m\ge 0\},
\\
\tau^+      =\, & \min\{m\in\mathbb{N}^*: X_m>0\},
\\
\tau^-      =\, & \min\{m\in\mathbb{N}^*: X_m<0\},
\\
\tau^{\pm}  =\, & \min\{m\in\mathbb{N}^*: X_m\neq0\},
\\
\tilde{\tau}^{\pm} =\, & \min\{m\ge \tau^{\rm o}: X_m\neq 0\}.
\end{align*}

We make the following assumptions on the sets $E^{\dag}$ and $E^{\rm o}$:
\begin{itemize}
\item[$(A_1)$]
if $X_0\in E^-$, then $\tau^{\rm o}=\tau^{\dag}$. This means that the chain
starting out of $E^{\dag}$ enters $E^{\dag}$ necessarily by passing through $E^{\rm o}$;
\item[$(A_2)$]
if $X_0\in E^+$, then $\tau^{\rm o}\le\tau^--1$ or,
equivalently, $\tau^-\ge\tau^{\rm o}+1$. This means that the chain starting in
$E^+$ exits $E^{\dag}$ necessarily by passing through $E^{\rm o}$.
\end{itemize}
Roughly speaking, $E^{\rm o}$ acts as a kind of `boundary' of $E^{\dag}$,
while $E^+$ acts as a kind of `interior' of $E^{\dag}$.
These assumptions are motivated by the example of integer-valued $(L,R)$-random
walks for which jumps are bounded from above by $R$ and from below by $-L$
($L,R\in\mathbb{N}^*$, $L$ for \textit{left}, $R$ for \textit{right},
the jumps lie in $\{-L,-L+1,\dots,R-1,R\}$).
If we consider the sojourn time above level~$0$, we are naturally
dealing with a `thick' boundary above $0$: $E^{\rm o}=\{0,1,\dots,M\}$ where $M$
is the maximum of $L$ and $R$. Section~\ref{section-LR} is devoted
to this class of random walks.

\subsection{Sojourn times}\label{subsection-sojourn-times}

We consider the sojourn time of $(X_m)_{m\in\mathbb{N}}$ in $E^{\dag}$ up to a fixed
time $n\in\mathbb{N}$: $T_0=0$ and, for $n\ge 1$,
$$
T_n=\#\{m\in\{1,\dots n\}: X_m\in E^{\dag}\}=\sum_{m=1}^n\mathbbm{1}_{E^{\dag}}(X_m).
$$
The random variable $T_n$ counts the indices $m\in\{1,\dots,n\}$
for which $X_m\in E^{\dag}$. Of course, we have $0\le T_n\le n$.
Inspired by the observation mentioned within the introduction
concerning the alternative counting process~(\ref{modified-sojourn}),
we define another sojourn time consisting in counting
the $m\in\{1,\dots,n\}$ such that
\\
\begin{itemize}
\item either $X_m\in E^+$,
\item or $X_m\in E^{\rm o}$ and $X_{m-1}\in E^+$,
\item or $X_m, X_{m-1}\in E^{\rm o}$, $X_{m-2}\in E^+$,
\item[$\vdots$]
\item or $X_{m},\dots, X_{2}\in E^{\rm o}$, $X_1\in E^+$.
\end{itemize}
Roughly speaking, we count the $X_m$ lying in $E^+$ or the $X_m$
lying in $E^{\rm o}$ `coming'  from a previous point lying in $E^+$.
Let us introduce the following events: $B_1=A_{0,1}=\{X_1\in E^+\}$
and for any $m\in\mathbb{N}\backslash\{0,1\}$,
\begin{align*}
A_{0,m}=\, & \{X_m\in E^+\},
\\
A_{\ell,m}=\, & \{X_m, X_{m-1},\dots, X_{m-\ell+1}\in E^{\rm o}\}
\cap\{X_{m-\ell}\in E^+\}\quad\text{for } \ell\in\{1,\dots,m-1\},
\\
B_m=\, & A_{0,m}\cup A_{1,m}\cup\dots\cup A_{m-1,m}.
\end{align*}
In other terms, we consider the sojourn time defined by $\tilde T_0=0$ and
for $n\in\mathbb{N}^*$, setting $\delta_m=\mathbbm{1}_{B_m}$,
$$
\tilde{T}_n=\#\{m\in\{1,\dots n\}: \delta_m=1\}=\sum_{m=1}^n\delta_m.
$$

It is interesting to know when, conversely, we do not count $X_m$ through this process.
This boils down to characterize the complement of $B_m$.
For this, we set also $B'_1=A'_{0,1}=\{X_1\in E^-\cup E^{\rm o}\}$
and for any $m\in\mathbb{N}\backslash\{0,1\}$,
\begin{align*}
A'_{0,m}=\, & \{X_m\in E^-\},
\\
A'_{\ell,m}=\, & \{X_m, X_{m-1},\dots, X_{m-\ell+1}\in E^{\rm o}\}
\cap\{X_{m-\ell}\in E^-\}\quad\text{for } \ell\in\{1,\dots,m-2\},
\\
A'_{m-1,m}=\, & \{X_m, X_{m-1},\dots, X_2\in E^{\rm o}\}\cap\{X_1\in E^-\cup E^{\rm o}\},
\\
B'_m=\, & A'_{0,m}\cup A'_{1,m}\cup\dots\cup A'_{m-1,m}.
\end{align*}
We have the following property the proof of which is postponed
to Section~\ref{proof-section}.
\begin{prop}\label{lemma}
For any $m\in\mathbb{N}^*$, the set $B'_m$ is the complementary of $B_m$: $B'_m=B_m^c$.
\end{prop}
\begin{rem}\label{remark-disjoint}
The sets $A_{\ell,m}$, $0\le\ell\le m-1$, are disjoint two by two, and the same
holds for $A'_{\ell,m}$, $0\le\ell\le m-1$. Hence, from Proposition~\ref{lemma},
we deduce the following identities:
$$
\delta_m=\sum_{\ell=0}^m\mathbbm{1}_{A_{\ell,m}}=\mathbbm{1}_{B_m}
=1-\mathbbm{1}_{B_m^c}=1-\mathbbm{1}_{B'_m}=1-\sum_{\ell=0}^m\mathbbm{1}_{A'_{\ell,m}}.
$$
\end{rem}

Our aim is to provide a methodology for deriving the joint probability distributions
of $(T_n,X_n)$ and  $(\tilde{T}_n,X_n)$.
For this, we develop a method for computing the generating functions
$\sum_{m,n\in\mathbb{N}:\atop m\le n}\mathbb{P}_i\{T_n=m, X_n\in F\}\,x^m y^{n-m}$
(Subsection~\ref{subsection-generating-Tn}) and
$\sum_{m,n\in\mathbb{N}:\atop m\le n}\mathbb{P}_i\{\tilde{T}_n=m, X_n\in F\}\,x^m y^{n-m}$
(Subsection~\ref{subsection-generating-Tn-tilde}) for any subset $F$ of $\mathcal{E}$.
When $F=\mathcal{E}$, these quantities are simply related to the probability distributions
of $T_n$ and $\tilde{T}_n$.
When $F=\{i\}$ for a fixed state $i\in\mathcal{E}$, this yields the probability
distributions of the sojourn times in $E^{\dag}$ up to time $n$ for the `bridge' (i.e.,
the pinned Markov chain) $(X_m|X_n=i)_{m\in\{0,\dots,n\}}$.

Actually, concerning $\tilde{T}_n$, we will only focus on the situation where
$X_0\in E^{\rm o},X_1\in E^{\pm}$ for lightening the paper and facilitating
the reading. That is, we will only provide a way for
computing the probabilities $\mathbb{P}_i\{X_1\in E^{\pm},\tilde{T}_n=m, X_n\in F\}$
for $i\in E^{\rm o}$ (Theorem~\ref{th-Ktildeixy}).
The study of $\tilde{T}_n$ subject to the complementary conditions $X_0,X_1\in E^{\rm o}$
or $X_0\notin E^{\rm o}$ could be carried out by using the results obtained
under the previous conditions. We explain below with few details how to proceed.
\begin{itemize}
\item
Under the condition $X_0,X_1\in E^{\rm o}$, we consider the first exit time from $E^{\rm o}$,
namely $\tau^{\pm}$. Then, the first terms of the chain satisfy
$X_0,X_1,\dots, X_{\tau^{\pm}-1}\in E^{\rm o}$ and
$X_{\tau^{\pm}}\in E^{\pm}$ and we have $\tilde{T}_{\tau^{\pm}-1}=0$.
Viewing the chain from time $\tau^{\pm}-1$, the facts $X_{\tau^{\pm}-1}\in E^{\rm o},
X_{\tau^{\pm}}\in E^{\pm}$ are the conditions analogous to $X_0\in E^{\rm o},X_1\in E^{\pm}$.

\item
Under the condition $X_0\notin E^{\rm o}$, we consider the first entrance time in $E^{\rm o}$,
and next the first exit time from $E^{\rm o}$, namely $\tilde{\tau}^{\pm}$. Then,
the first terms of the chain satisfy
$X_0,X_1,\dots, X_{\tau^0-1}\in E^{\pm}$, $X_{\tau^0},\dots,
X_{\tilde{\tau}^{\pm}-1}\in E^{\rm o}$ and $X_{\tilde{\tau}^{\pm}}\in E^{\pm}$
and we have $\tilde{T}_{\tilde{\tau}^{\pm}-1}=\tilde{\tau}^{\pm}-1$
if $X_0\in E^+$ and $\tilde{T}_{\tilde{\tau}^{\pm}-1}=0$ if $X_0\in E^-$.
From time $\tilde{\tau}^{\pm}-1$, the facts $X_{\tilde{\tau}^{\pm}-1}\in E^{\rm o},
X_{\tilde{\tau}^{\pm}}\in E^{\pm}$ are the conditions analogous to $X_0\in E^{\rm o},X_1\in E^{\pm}$.
\end{itemize}
In both situations, we could perform the computations by appealing to
identities similar to~(\ref{P-tildeTn-m}) (which are used for proving
Theorem~\ref{th-Ktildeixy}) with the help of Markov property.
We will let the reader write the details.

\section{Generating functions}\label{section-gene-func}

\subsection{Generating functions of $X,\tau^{\rm o},\tau^{\dag},\tau^+,\tau^-$}

Let us introduce the following generating functions: for any states
$i,j\in\mathcal{E}$, any subset $F\subset\mathcal{E}$ and any real number $x$
such that the following series converge,
\\
\begin{align*}
G_{ij}(x)=\;
& \sum_{m=0}^{\infty}\mathbb{P}_i\{X_m=j\}\,x^m,\quad
G_i(x)=\sum_{m=0}^{\infty}\mathbb{P}_i\{X_m\in F\}\,x^m=\sum_{j\in F} G_{ij}(x),
\\
H_{ij}^{\rm o}(x)=\;
& \sum_{m=1}^{\infty}\mathbb{P}_i\{\tau^{\rm o}=m,
X_{\tau^{\rm o}}=j\}\,x^m
=\mathbb{E}_i\big(x^{\tau^{\rm o}}\mathbbm{1}_{\{X_{\tau^{\rm o}}=j\}}\big),
\\
H_{ij}^{\dag}(x)=\;
& \mathbb{E}_i\big(x^{\tau^{\dag}}\mathbbm{1}_{\{X_{\tau^{\dag}}=j\}}\big),\quad
H_{ij}^+(x)=\mathbb{E}_i \big(x^{\tau^+}\mathbbm{1}_{\{X_{\tau^+}=j\}}\big),\quad
H_{ij}^-(x)=\mathbb{E}_i \big(x^{\tau^-}\mathbbm{1}_{\{X_{\tau^-}=j\}}\big).
\end{align*}
In the above notations, the indices $i$ and $j$ refer to starting points and
arrival points or entrance locations while superscripts $\rm o,\dag,+,-$ refer to
the entrance in the respective sets $E^{\rm o},E^{\dag},E^+,E^-$.
Moreover, for lightening the settings, when writing $X_{\tau^u}=j$ in $H_{ij}^u(x)$,
$u\in\{\rm o,\dag,+,-\}$, we implicitly restrict this event to the condition that $\tau^u<+\infty$
and we omit to write this condition explicitly.

For determining $G_{ij}(x)$, the standard method consists in using the
well-known Chapman-Kolmogorov equation:
for any $m\in\mathbb{N}$, $\ell\in\{0,\dots,m\}$ and $i,j\in\mathcal{E}$,
\begin{equation}\label{CK}
\mathbb{P}_i\{X_m=j\}=\sum_{k\in\mathcal{E}}
\mathbb{P}_i\{X_{\ell}=k\}\,\mathbb{P}_k\{X_{m-\ell}=j\}.
\end{equation}
In particular, by choosing $\ell=1$ in~(\ref{CK}), we have, for $m\ge 1$,
that
$$
\mathbb{P}_i\{X_m=j\}=\sum_{k\in\mathcal{E}} p_{ik}
\,\mathbb{P}_k\{X_{m-1}=j\}
$$
and we get that
$$
G_{ij}(x)=\delta_{ij}+\sum_{m=1}^{\infty} \!\Bigg(\sum_{k\in\mathcal{E}}
p_{ik}\,\mathbb{P}_k\{X_{m-1}=j\}\Bigg) x^m
=\delta_{ij}+\sum_{k\in\mathcal{E}} p_{ik}
\sum_{m=1}^{\infty}\mathbb{P}_k\{X_{m-1}=j\}\,x^m.
$$
We obtain the famous backward Kolmogorov equation
\begin{equation}\label{eqG1}
G_{ij}(x)=\delta_{ij}+ x\sum_{k\in\mathcal{E}} p_{ik} G_{kj}(x)
\quad\text{for } i,j\in\mathcal{E}.
\end{equation}
Similarly, by choosing $\ell=m-1$ with $m\ge 1$ in~(\ref{CK}), we get the
forward Kolmogorov equation
\begin{equation}\label{eqG2}
G_{ij}(x)=\delta_{ij}+x\sum_{k\in\mathcal{E}} p_{kj} G_{ik}(x)
\quad\text{for } i,j\in\mathcal{E}.
\end{equation}
For determining $H_{ij}^{\rm o}(x)$, we observe that
for $i\in\mathcal{E}$, $j\in E^{\rm o}$ and $m\in\mathbb{N}^*$,
if $X_0=i$ and $X_m=j$, then the chain has entered $E^{\rm o}$ between
times $1$ and $m$. In symbols, appealing to the strong Markov property,
$$
\mathbb{P}_i\{X_m=j\}=\mathbb{P}_i\{\tau^{\rm o}\le m, X_m=j\}
=\sum_{\ell=1}^m\sum_{k\in E^{\rm o}}\mathbb{P}_i\{\tau^{\rm o}=\ell,
X_{\tau^{\rm o}}=k\}\,\mathbb{P}_k\{X_{m-\ell}=j\}.
$$
Therefore,
\begin{align*}
G_{ij}(x)=\;
& \mathbb{P}_i\{X_0=j\}+
\sum_{m=1}^{\infty} \!\Bigg(\sum_{\ell=1}^m\sum_{k\in E^{\rm o}}
\mathbb{P}_i\{\tau^{\rm o}=\ell, X_{\tau^{\rm o}}=k\}
\,\mathbb{P}_k\{X_{m-\ell}=j\}\Bigg)x^m
\\
=\; & \delta_{ij}+\sum_{\ell=1}^{\infty}\sum_{k\in E^{\rm o}}
\mathbb{P}_i\{\tau^{\rm o}=\ell, X_{\tau^{\rm o}}=k\}\sum_{m=\ell}^{\infty}
\mathbb{P}_k\{X_{m-\ell}=j\}\,x^m
\\
=\; & \delta_{ij}+\sum_{\ell=1}^{\infty}\sum_{k\in E^{\rm o}}
\mathbb{P}_i\{\tau^{\rm o}=\ell, X_{\tau^{\rm o}}=k\}\,x^{\ell} G_{kj}(x)
\\
=\; & \delta_{ij}+\sum_{k\in E^{\rm o}}\mathbb{E}_i \big(x^{\tau^{\rm o}}
\mathbbm{1}_{\{ X_{\tau^{\rm o}}=k\}}\big) G_{kj}(x).
\end{align*}
We get the equation
\begin{equation}\label{eqGH1}
G_{ij}(x)=\delta_{ij}+\sum_{k\in E^{\rm o}} H_{ik}^{\rm o}(x) G_{kj}(x)
\quad\text{for } i\in\mathcal{E},j\in E^{\rm o}.
\end{equation}
Let us point out that, due to Assumptions~$(A_1)$ and $(A_2)$, (\ref{eqGH1})
holds true for $i\in E^+$, $j\in E^-$, and also
for $i\in E^-$, $j\in E^{\dag}$. In the same way, we have that
\begin{align}
G_{ij}(x)=\;
& \delta_{ij}+\sum_{k\in E^{\dag}} H_{ik}^{\dag}(x) G_{kj}(x)
\quad\text{for } i\in\mathcal{E},j\in E^{\dag},
\label{eqGH2dag}\\
G_{ij}(x)=\;
& \delta_{ij}+\sum_{k\in E^+} H_{ik}^+(x) G_{kj}(x)
\quad\text{for } i\in\mathcal{E},j\in E^+,
\label{eqGH2+}\\
G_{ij}(x)=\;
& \delta_{ij}+\sum_{k\in E^-} H_{ik}^-(x) G_{kj}(x)
\quad\text{for } i\in\mathcal{E},j\in E^-.
\label{eqGH2-}
\end{align}

Moreover, if $i\in E^-$, then, by Assumption~$(A_1)$, $\tau^{\dag}=\tau^{\rm o}$ which
entails that $H_{ij}^{\dag}(x)=H_{ij}^{\rm o}(x)$.
If $i\in E^+$, then, by Assumption~$(A_2)$, $\tau^{\dag}=1$ which
entails that $H_{ij}^{\dag}(x)=p_{ij}\,x$. If $i\in E^{\rm o}$,
$$
H_{ij}^{\dag}(x)=\sum_{k\in E^{\dag}}\mathbb{E}_i \big(x^{\tau^{\dag}}\mathbbm{1}_{\{X_1=k,
X_{\tau^{\dag}}=j\}}\big)+\sum_{k\in E^-}\mathbb{E}_i \big(x^{\tau^{\dag}}
\mathbbm{1}_{\{X_1=k,X_{\tau^{\dag}}=j\}}\big).
$$
If $X_1\in E^{\dag}$, then $\tau^{\dag}=1$ while if $X_1\in E^-$, then
$\tau^{\dag}=\tau^{\rm o}$.
Now, let us introduce the first hitting time of $E^{\dag}$ by
$(X_m)_{m\in\mathbb{N}}$: $\tau^{\dag\prime}=\min\{m\in\mathbb{N}: X_m\in E^{\dag}\}$.
Of course, we have $\tau^{\dag\prime}=0$ if $X_0\in E^{\dag}$,
$\tau^{\dag\prime}=\tau^{\dag}$ if $X_0\in E^-$ and $\tau^{\dag\prime}$
is related to $\tau^{\dag}$ according to $\tau^{\dag}=1+\tau^{\dag\prime}\circ\theta_1$
where $\theta_1$ is the usual shift operator (acting as $X_m\circ\theta_1=X_{m+1}$).
Moreover, $X_{\tau^{\dag}}=X_{\tau^{\dag\prime}}\circ\theta_1$.
With these settings at hands, thanks to the Markov property, we obtain
$\mathbb{E}_i \big(x^{\tau^{\dag}}\mathbbm{1}_{\{X_1=k,X_{\tau^{\dag}}=j\}}\big)
= p_{ik}\,x\,\mathbb{E}_k \big(x^{\tau^{\dag\prime}}
\mathbbm{1}_{\{ X_{\tau^{\dag\prime}}=j\}}\big)$ which yields that
$$
\mathbb{E}_i \big(x^{\tau^{\dag}}\mathbbm{1}_{\{X_1=k,X_{\tau^{\dag}}=j\}}\big)
=\begin{cases}
\delta_{jk}\,p_{ik}\,x & \text{if $k\in E^{\dag}$,}
\\
p_{ik}\,x \,\mathbb{E}_k \big(x^{\tau^{\rm o}}\mathbbm{1}_{\{X_{\tau^{\rm o}}=j\}}\big)
& \text{if $k\in E^-$.}
\end{cases}
$$
As a result, we see that the function $H_{ij}^{\dag}$ can be expressed by means of
$H_{ij}^{\rm o}$ according to
\begin{equation}\label{link-H-H+}
H_{ij}^{\dag}(x)=\begin{cases}
H_{ij}^{\rm o}(x) & \text{if $i\in E^-$,}
\\
p_{ij}\,x & \text{if $i\in E^+$,}
\\[1ex]
\displaystyle x\Bigg(p_{ij}+\sum_{k\in E^-} p_{ik} H_{kj}^{\rm o}(x)\Bigg)
& \text{if $i\in E^{\rm o}$.}
\end{cases}
\end{equation}

In the sequel of the paper, we will also use the generating functions below.
Set, for any $i\in\mathcal{E}$ and $j\in E^{\rm o}$,
\begin{align*}
H_{ij}^{\rm o\dag}(x)=\;
& \mathbb{E}_i \big(x^{\tau^{\rm o}}
\mathbbm{1}_{\{X_1\in E^{\dag},X_{\tau^{\rm o}}=j\}}\big),
\\
H_{ij}^{\rm o+}(x)=\;
& \mathbb{E}_i \big(x^{\tau^{\rm o}}
\mathbbm{1}_{\{X_1\in E^+,X_{\tau^{\rm o}}=j\}}\big),
\\
H_{ij}^{\rm o-}(x)=\;
& \mathbb{E}_i\big(x^{\tau^{\rm o}}
\mathbbm{1}_{\{X_1\in E^-,X_{\tau^{\rm o}}=j\}}\big).
\end{align*}
In short, $H_{ij}^{uv}(x)$, $u,v\in\{\rm o,\dag,+,-\}$, is related to time $\tau^u$
and the first step $X_1\in E^v$.
We propose a method for calculating $H_{ij}^{\rm o\dag}(x)$,
$H_{ij}^{\rm o+}(x)$ and $H_{ij}^{\rm o-}(x)$.
For this, we introduce the first hitting time of $E^{\rm o}$ by
$(X_m)_{m\in\mathbb{N}}$: $\tau^{\rm o\prime}=\min\{m\in\mathbb{N}: X_m\in E^{\rm o}\}$.
Of course, we have $\tau^{\rm o\prime}=0$ if $X_0\in E^{\rm o}$,
$\tau^{\rm o\prime}=\tau^{\rm o}$ if $X_0\notin E^{\rm o}$ and $\tau^{\rm o\prime}$
is related to $\tau^{\rm o}$ according to $\tau^{\rm o}=1+\tau^{\rm o\prime}\circ\theta_1$.
Moreover, $X_{\tau^{\rm o}}=X_{\tau^{\rm o\prime}}\circ\theta_1$.
With these settings at hands, thanks to the Markov property, we obtain,
for any $i\in\mathcal{E}$ and $j\in E^{\rm o}$, that
$$
H_{ij}^{\rm o\dag}(x)=x\sum_{k\in E^{\dag}} p_{ik}\,\mathbb{E}_k
\big(x^{\tau^{\rm o\prime}}\mathbbm{1}_{\{ X_{\tau^{\rm o\prime}}=j\}}\big)
=x\Bigg(\sum_{k\in E^{\rm o}} p_{ik}\delta_{kj}
+\sum_{k\in E^+} p_{ik} \,\mathbb{E}_k
\big(x^{{\tau^{\rm o}}}\mathbbm{1}_{\{ X_{{\tau^{\rm o}}}=j\}}\big)\Bigg)
$$
which simplifies into
\begin{equation}\label{Ho-dag}
H_{ij}^{\rm o\dag}(x)=x\Bigg(p_{ij}+\sum_{k\in E^+} p_{ik} H_{kj}^{\rm o}(x)\Bigg)\!.
\end{equation}
Similarly,
\begin{equation}\label{Ho+-}
H_{ij}^{\rm o+}(x)=x\sum_{k\in E^+} p_{ik} H_{kj}^{\rm o}(x),\quad
H_{ij}^{\rm o-}(x)=x\sum_{k\in E^-} p_{ik} H_{kj}^{\rm o}(x).
\end{equation}

\subsection{Generating function of $T_n$}\label{subsection-generating-Tn}

Now, we introduce the generating function of the numbers
$\mathbb{P}_i\{T_n=m, X_n\in F\}$, $m,n\in\mathbb{N}$: set, for any
$i\in\mathcal{E}$ and any real numbers $x,y$ such that the following series
converges,
$$
K_i(x,y)=\sum_{m,n\in\mathbb{N}:\atop m\le n}
\mathbb{P}_i\{T_n=m, X_n\in F\}\,x^m y^{n-m}.
$$
In the theorem below, we provide a way for computing $K_i(x,y)$.
%
\begin{teo}\label{th-Kixy}
The generating function $K_i$, $i\in\mathcal{E}$, satisfies the equation
\begin{equation}\label{eqKxy}
K_i(x,y)=K_i(x,0)+K_i(0,y)+\sum_{j\in E^{\rm o}} \!\left(H_{ij}^{\rm o\dag}(x)
+\frac x y\,H_{ij}^{\rm o-}(y)\right)\! K_j(x,y)-\sum_{j\in E^{\rm o}} H_{ij}^{\rm o\dag}(x) K_j(x,0)
-\mathbbm{1}_F(i)
\end{equation}
where, for any $i\in\mathcal{E}$,
\begin{equation}\label{eqK0xy}
K_i(x,0)=G_i(x)-\sum_{j\in E^-} H_{ij}^-(x) G_j(x),
\quad K_i(0,y)=G_i(y)-\sum_{j\in E^{\dag}} H_{ij}^{\dag}(y) G_j(y),
\end{equation}
and where the functions
$H_{ij}^-$, $j\in E^-$, and $H_{ij}^{\dag}$, $j\in E^{\dag}$, are
given by~(\ref{eqGH2dag}) and~(\ref{eqGH2-}), and the functions
$H_{ij}^{\rm o\dag}$ and $H_{ij}^{\rm o-}$, $j\in E^{\rm o}$ are given
by~(\ref{Ho-dag}) and~(\ref{Ho+-}).
\end{teo}
%

%
\begin{rem}\label{F=space1}
If we choose $F=\mathcal{E}$, for all state $i$ in $\mathcal{E}$,
we simply have $G_i(x)=\frac{1}{1-x}$ and~(\ref{eqK0xy}) yields that
$$
K_i(x,0)=\frac{1-\mathbb{E}_i\big(x^{\tau^-}\big)}{1-x},\quad
K_i(0,y)=\frac{1-\mathbb{E}_i\big(y^{\tau^{\dag}}\big)}{1-y}.
$$
\end{rem}
%
\begin{rem}\label{remark-one-point1}
If $E^{\rm o}$ reduces to one point $i_0$, then, for $i=i_0$, (\ref{eqKxy}) immediately
yields
$$
K_{i_0}(x,y)=\frac{\big(1-H_{i_0i_0}^{\rm o\dag}(x)\big)K_{i_0}(x,0)+K_{i_0}(0,y)
-\mathbbm{1}_F(i_0)}{1-H_{i_0i_0}^{\rm o\dag}(x)-\frac x y\,H_{i_0i_0}^{\rm o-}(y)}.
$$
Moreover, $H_{i_0i_0}^{\rm o\dag}(x)$ and $H_{i_0i_0}^{\rm o-}(y)$ can be computed thanks
to~(\ref{Ho-dag}) and (\ref{Ho+-}) with the aid of $H_{ki_0}^{\rm o}(x)=G_{ki_0}(x)/G_{i_0i_0}(x)$ which
comes from~(\ref{eqGH1}).
\end{rem}
%
Let us mention that due to~(\ref{eqKxy}), it is enough to know $K_i(x,y)$
only for $i\in E^{\rm o}$ to derive $K_i(x,y)$ for $i\in\mathcal{E}\backslash E^{\rm o}$.
Indeed, we have the connections below.
%
\begin{prop} \label{propKixy}
For $i\in\mathcal{E}\backslash E^{\rm o}$, $K_i(x,y)$ can be expressed
by means of $K_j(x,y)$, $j\in E^{\rm o}$, according to the following
identities: if $i\in E^+$,
\begin{equation}\label{uno}
K_i(x,y)=G_i(x)-\sum_{j\in E^{\rm o}} H_{ij}^{\rm o}(x) G_j(x)
+\sum_{j\in E^{\rm o}} H_{ij}^{\rm o}(x) K_{j}(x,y)
\end{equation}
and if $i\in E^-$,
\begin{equation}\label{due}
K_i(x,y)=G_i(y)-\sum_{j\in E^{\rm o}} H_{ij}^{\rm o}(y) G_j(y)
+\frac x y\sum_{j\in E^{\rm o}} H_{ij}^{\rm o}(y) K_j(x,y).
\end{equation}
\end{prop}

\subsection{Generating functions of $\tau^{\pm}$ and $\tilde{\tau}^{\pm}$}
\label{subsection-gene}

From now on we restrict ourselves to the case where $X_0\in E^{\rm o}$ and $X_1\in E^{\pm}$.
In Subsection~\ref{subsection-sojourn-times}, we indicate how to treat the
complementary case. Let us introduce the following generating functions: for any $i,j\in E^{\rm o}$,
\begin{align*}
H_{ij}^{\pm}(x)=\;
& \mathbb{E}_i \big(x^{\tau^{\pm}-1}\mathbbm{1}_{\{X_{\tau^{\pm}-1}=j\}}\big)
=\sum_{m=1}^{\infty} \mathbb{P}_i\{\tau^{\pm}=m,X_{\tau^{\pm}-1}=j\}\,x^{m-1},
\\
\tilde{H}_{ij}^{\pm+}(x)
=\; & \mathbb{E}_i \big(x^{\tilde{\tau}^{\pm}-1}
\mathbbm{1}_{\{X_1\in E^+,X_{\tilde{\tau}^{\pm}-1}=j\}}\big)
=\sum_{m=1}^{\infty} \mathbb{P}_i\{X_1\in E^+,
\tilde{\tau}^{\pm}=m,X_{\tilde{\tau}^{\pm}-1}=j\}\,x^{m-1},
\\
\tilde{H}_{ij}^{\pm-}(x)=\;
& \mathbb{E}_i \big(x^{\tilde{\tau}^{\pm}-1}
\mathbbm{1}_{\{X_1\in E^-,X_{\tilde{\tau}^{\pm}-1}=j\}}\big)
=\sum_{m=1}^{\infty} \mathbb{P}_i\{X_1\in E^-,
\tilde{\tau}^{\pm}=m,X_{\tilde{\tau}^{\pm}-1}=j\}\,x^{m-1}.
\end{align*}
By noticing that $\tilde{\tau}^{\pm}=\tau^{\rm o}+\tau^{\pm}\circ\theta_{\tau^{\rm o}}$
where $\theta_{\tau^{\rm o}}$ acts on $(X_m)_{m\in\mathbb{N}}$
as $X_m\circ \theta_{\tau^{\rm o}}=X_{m+\tau^{\rm o}}$
and by using the strong Markov property, we have that
$$
\tilde{H}_{ij}^{\pm+}(x)=\sum_{k\in E^{\rm o}} \mathbb{E}_i \big(x^{\tau^{\rm o}}
\mathbbm{1}_{\{X_1\in E^+,X_{\tau^{\rm o}}=k\}}\big)\,
\mathbb{E}_k \big(x^{\tau^{\pm}-1}\mathbbm{1}_{\{X_{\tau^{\pm}-1}=j\}}\big),
$$
and a similar expression holds for $\tilde{H}_{ij}^{\pm-}(x)$. In short, we have obtained that
$$
\tilde{H}_{ij}^{\pm+}(x)=\sum_{k\in E^{\rm o}} H_{ik}^{\rm o+}(x) H_{kj}^{\pm}(x),\quad
\tilde{H}_{ij}^{\pm-}(x)=\sum_{k\in E^{\rm o}} H_{ik}^{\rm o-}(x) H_{kj}^{\pm}(x),
$$
and, by~(\ref{Ho+-}),
\begin{equation}\label{tilde-K*}
\tilde{H}_{ij}^{\pm+}(x)=x\sum_{k\in E^+,\ell\in E^{\rm o}} p_{ik}
H_{k\ell}^{\rm o}(x)H_{\ell j}^{\pm}(x),\quad
\tilde{H}_{ij}^{\pm-}(x)=x\sum_{k\in E^-,\ell\in E^{\rm o}} p_{ik}
H_{k\ell}^{\rm o}(x)H_{\ell j}^{\pm}(x).
\end{equation}

Hence, we need to evaluate $H_{ij}^{\pm}(x)$ for $i,j\in E^{\rm o}$.
We note that, if $X_1\in E^{\pm}$, then $\tau^{\pm}=1$
and $X_{\tau^{\pm}-1}=X_0$, so that we get
\begin{align*}
H_{ij}^{\pm}(x)=\;
& \mathbb{E}_i \big(x^{\tau^{\pm}-1}\mathbbm{1}_{\{X_1\in E^{\pm},
X_{\tau^{\pm}-1}=j\}}\big)+\mathbb{E}_i \big(x^{\tau^{\pm}-1}
\mathbbm{1}_{\{X_1\in E^{\rm o},X_{\tau^{\pm}-1}=j\}}\big)
\\
=\; & \delta_{ij}\mathbb{P}_i\{X_1\in E^{\pm}\}
+x\sum_{k\in E^{\rm o}} \mathbb{P}_i \{X_1=k\}\,
\mathbb{E}_k \big(x^{\tau^{\pm}-1}\mathbbm{1}_{\{X_{\tau^{\pm}-1}=j\}}\big).
\end{align*}
As a result, we obtain the following equation: for $i,j\in E^{\rm o}$,
\begin{equation}\label{eqHtildedag-}
H_{ij}^{\pm}(x)=\delta_{ij}\mathbb{P}_i\{X_1\in E^{\pm}\}
+x\sum_{k\in E^{\rm o}}p_{ik} H_{kj}^{\pm}(x).
\end{equation}
%
\begin{rem}\label{remark-one-point2}
If $E^{\rm o}$ reduces to one point $i_0$, then (\ref{eqHtildedag-})
immediately yields $H_{i_0i_0}^{\pm}(x)=(1-p_{i_0i_0})/(1-p_{i_0i_0}x)$
which in turn yields that
$$
\tilde{H}_{i_0i_0}^{\pm+}(x)=\frac{1-p_{i_0i_0}}{1-p_{i_0i_0}x}\,H_{i_0i_0}^{\rm o+}(x),\quad
\tilde{H}_{i_0i_0}^{\pm-}(x)=\frac{1-p_{i_0i_0}}{1-p_{i_0i_0}x}\,H_{i_0i_0}^{\rm o-}(x).
$$
If we additionally impose the (more restrictive) condition that the Markov chain does not stay
at its current location in $E^{\rm o}$, that is, $p_{ii}=0$ for any $i\in E^{\rm o}$,
then we simply have $H_{i_0i_0}^{\pm}(x)=1$,
$\tilde{H}_{i_0i_0}^{\pm+}(x)=H_{i_0i_0}^{\rm o+}(x)$ and
$\tilde{H}_{i_0i_0}^{\pm-}(x)=H_{i_0i_0}^{\rm o-}(x)$.
\end{rem}

\subsection{Generating functions of $X,\tau^+,\tau^-$ subjected to $X_1\in E^{\pm}$}

In what follows we need the generating functions below: for any $i\in E^{\rm o}$, we set
\begin{align*}
G_i^+(x)=\; & \sum_{m=0}^{\infty} \mathbb{P}_i\{X_1\in E^+,X_m\in F\}\,x^m,
\\
G_i^-(x)=\; & \sum_{m=0}^{\infty} \mathbb{P}_i\{X_1\in E^-,X_m\in F\}\,x^m;
\\
\noalign{\noindent for $i\in E^{\rm o},j\in E^-$,}
H_{ij}^{-+}(x)&=\sum_{m=1}^{\infty} \mathbb{P}_i\{X_1\in E^+,
\tau^-=m,X_{\tau^-}=j\}\,x^m=\mathbb{E}_i \big(x^{\tau^-}
\mathbbm{1}_{\{X_1\in E^+,X_{\tau^-}=j\}}\big);
\\
\noalign{\noindent for $i\in E^{\rm o},j\in E^+$,}
H_{ij}^{+-}(x)=\; & \sum_{m=1}^{\infty} \mathbb{P}_i\{X_1\in E^-,
\tau^{+}=m,X_{\tau^{+}}=j\}\,x^m=\mathbb{E}_i \big(x^{\tau^{+}}
\mathbbm{1}_{\{X_1\in E^-,X_{\tau^{+}}=j\}}\big).
\end{align*}
By the Markov property, we clearly have that
\begin{equation}\label{G+-}
G_i^+(x)=\mathbbm{1}_F(i) \mathbb{P}_i\{X_1\in E^+\}+x\sum_{j\in E^+} p_{ij} G_j(x),\quad
G_i^-(x)=\mathbbm{1}_F(i) \mathbb{P}_i\{X_1\in E^-\}+x\sum_{j\in E^-} p_{ij} G_j(x)
\end{equation}
and
\begin{equation}\label{H+-}
H_{ij}^{-+}(x)=x\sum_{k\in E^+} p_{ik} H_{kj}^-(x),\quad
H_{ij}^{+-}(x)=x\sum_{k\in E^-} p_{ik} H_{kj}^+(x).
\end{equation}

\subsection{Generating function of $\tilde{T}_n$}\label{subsection-generating-Tn-tilde}

Now, we introduce the generating function of the numbers
$\mathbb{P}_i\{X_1\in E^{\pm}, \tilde{T}_n=m, X_n\in F\}$, $m,n\in\mathbb{N}$:
for any $i\in E^{\rm o}$,
$$
\tilde{K}_i(x,y)=\sum_{m,n\in\mathbb{N}:\atop m\le n}
\mathbb{P}_i\{X_1\in E^{\pm},\tilde{T}_n=m, X_n\in F\}\,x^m y^{n-m}.
$$
In the following theorem, we provide a way for computing $\tilde{K}_i(x,y)$.
%
\begin{teo}\label{th-Ktildeixy}
The generating function $\tilde{K}_i$, $i\in E^{\rm o}$, satisfies the equation
\begin{align}
\tilde{K}_i(x,y)=\;
& \tilde{K}_i(x,0)+\tilde{K}_i(0,y)+\sum_{j\in E^{\rm o}}
\big(\tilde{H}_{ij}^{\pm+}(x)+\tilde{H}_{ij}^{\pm-}(y)\big)\tilde{K}_j(x,y)
\nonumber\\
& -\sum_{j\in E^{\rm o}}\tilde{H}_{ij}^{\pm+}(x)\tilde{K}_j(x,0)
-\sum_{j\in E^{\rm o}}\tilde{H}_{ij}^{\pm-}(y)\tilde{K}_j(0,y)
-\mathbbm{1}_F(i)\mathbb{P}_i\{X_1\in E^{\pm}\},
\label{eqKxy-tilde}
\end{align}
where, for any $i\in E^{\rm o}$,
\begin{align}
\tilde{K}_i(x,0)=\;
& G_i^+(x)-\sum_{j\in E^-} H_{ij}^{-+} (x) G_j(x)+\mathbbm{1}_F(i)\mathbb{P}_i\{X_1\in E^-\},
\nonumber\\[-3ex]
&
\label{tilde-K}
\\[-1ex]
\tilde{K}_i(0,y)=\;
& G_i^-(y)-\sum_{j\in E^+} H_{ij}^{+-}(y) G_j(y)+\mathbbm{1}_F(i)\mathbb{P}_i\{X_1\in E^+\},
\nonumber
\end{align}
and where the functions $\tilde{H}_{ij}^{\pm+}$ and $\tilde{H}_{ij}^{\pm-}$, $i,j\in E^{\rm o}$,
are given by~(\ref{tilde-K*}) and the functions $H_{ij}^{-+}$, $i\in E^{\rm o}$, $j\in E^-$, and
$H_{ij}^{+-}$, $i\in E^{\rm o}$, $j\in E^+$, are given by~(\ref{H+-}).
\end{teo}
%
\begin{rem}\label{F=space2}
If we choose $F=\mathcal{E}$, for all state $i$ in $\mathcal{E}$,
we simply have $G_i^+(x)=\mathbb{P}_i\{X_1\in E^+\}/(1-x)$
and $G_i^-(x)=\mathbb{P}_i\{X_1\in E^-\}/(1-x)$, and~(\ref{tilde-K}) yields that
\begin{align*}
\tilde{K}_i(x,0) =\,
& \frac{1}{1-x}\Big(\mathbb{P}_i\{X_1\in E^{\pm}\}-x\,\mathbb{P}_i\{X_1\in E^-\}
-x\sum_{j\in E^+} p_{ij}\,\mathbb{E}_j\big(x^{\tau^-}\big)\Big),
\\
\tilde{K}_i(0,y) =\,
& \frac{1}{1-y}\Big(\mathbb{P}_i\{X_1\in E^{\pm}\}-y\,\mathbb{P}_i\{X_1\in E^+\}
-y\sum_{j\in E^-} p_{ij}\,\mathbb{E}_j\big(y^{\tau^+}\big)\Big).
\end{align*}
\end{rem}
%

In view of Remark~\ref{remark-one-point2}, if $E^{\rm o}$ reduces to one point $i_0$,
then (\ref{eqKxy-tilde}) immediately yields the explicit expression below.
%
\begin{cor}\label{remark-one-point3}
If $E^{\rm o}$ reduces to one point $i_0$, then
$$
\tilde{K}_{i_0}(x,y)=\frac{\big(1-\tilde{H}_{i_0i_0}^{\pm+}(x)\big)\tilde{K}_{i_0}(x,0)
+\big(1-\tilde{H}_{i_0i_0}^{\pm-}(y)\big)\tilde{K}_{i_0}(0,y)-\mathbbm{1}_F(i_0)(1-p_{i_0i_0})}%
{1-\tilde{H}_{i_0i_0}^{\pm+}(x)-\tilde{H}_{i_0i_0}^{\pm-}(y)}.
$$
Moreover, under the additional assumption that the Markov chain does not
stay at its current location in $E^{\rm o}$, that is $p_{ii}=0$ for any
$i\in E^{\rm o}$, the generating function $\tilde{K}_{i_0}$
can be simplified into
$$
\tilde{K}_{i_0}(x,y)=\frac{\big(1-H^{\rm o+}_{i_0i_0}(x)\big)\tilde{K}_{i_0}(x,0)
+\big(1-H^{\rm o-}_{i_0i_0}(y)\big)\tilde{K}_{i_0}(0,y)-\mathbbm{1}_F(i_0)}%
{1-H^{\rm o+}_{i_0i_0}(x)-H^{\rm o-}_{i_0i_0}(y)}.
$$
\end{cor}
%

\section{Application to $(L,R)$-random walk}\label{section-LR}

The algebraic equations~(\ref{eqKxy}) and (\ref{eqKxy-tilde})
satisfied by the generating functions $K_i$ and $\tilde{K}_i$,
$i\in E^{\rm o}$, produce linear systems which may be rewritten
in a matrix form possibly involving infinite matrices.
In this section, we focus on the case of the $(L,R)$-random walk. In this case,
the systems of interest consist of a finite number of equations and may be solved
by using a matrix approach that we describe here.
We guess that this approach provides a methodology which should be efficiently
numerically implemented.

\subsection{Settings}

Let $L,R$ be positive integers, $M=\max(L,R)$ and let $(U_{\ell})_{\ell\in\mathbb{N}^*}$ be a
sequence of independent identically distributed random variables with values
in $\{-L,-L+1,\dots,R-1,R\}$. Put $\pi_i=\mathbb{P}\{U_1=i\}$ for
$i\in\{-L,\dots,R\}$ and $\pi_i=0$ for $i\in\mathbb{Z}\backslash\{-L,\dots,R\}$.
The common generating function of the $U_{\ell}$'s is given by
$$
\mathbb{E}\!\left(y^{U_1}\right)=\sum_{j=-L}^R \pi_j\,y^j
=y^{-L}\sum_{j=0}^{L+R} \pi_{j-L}\,y^j.
$$
Let $X_0$ be an integer and set $X_m=X_0+\sum_{\ell=1}^m U_{\ell}$ for any
$m\in\mathbb{N}^*$. Set $U_0=X_0$, and notice that $X_m$ is the partial sum
of the series $\sum U_{\ell}$. The Markov chain $(X_m)_{m\in\mathbb{N}}$ is a random walk
defined on $\mathcal{E}=\mathbb{Z}$ with transition probabilities defined
as $p_{ij}=\mathbb{P}\{X_{m+1}=j\,|\,X_m=i\}=\pi_{j-i}$. The jumps are bounded
and we have that
\begin{equation}\label{cancel}
p_{ij}=0 \text{ for } j\notin\{i-L,i-L+1,\dots, i+R\}.
\end{equation}
We choose here
\begin{align*}
E^{\rm o} =\, & \{0,1,\dots,M-1\},\\
E^{\dag}  =\, & \{0,1,\dots\},\\
E^+       =\, & \{M,M+1,\dots\},\\
E^-       =\, & \{\dots,-2,-1\}.
\end{align*}
The settings of Section~\ref{section-settings} can be rewritten in this context as
\begin{align*}
T_n         =\, & \#\{m\in\{1,\dots,n\}: X_m\ge 0\},
\\
\tau^{\rm o}=\, & \min\{m\in\mathbb{N}^*: X_m\in\{0,1,\dots,M-1\}\},
\\
\tau^{\dag} =\, & \min\{m\in\mathbb{N}^*: X_m\ge 0\},
\\
\tau^+      =\, & \min\{m\in\mathbb{N}^*: X_m\ge M\},
\\
\tau^-      =\, & \min\{m\in\mathbb{N}^*: X_m\le -1\}.
\end{align*}
It is easy to see that Assumptions~$(A_1)$ and $(A_2)$ are fulfilled. Indeed,
\begin{itemize}
\item
if $X_0<0$, then $X_0,X_1,\dots,X_{\tau^{\dag}-1}<0$ and $X_{\tau^{\dag}}\ge0$.
Since $X_{\tau^{\dag}}=X_{\tau^{\dag}-1}+U_{\tau^{\dag}}$, $X_{\tau^{\dag}-1}\le -1$ and
$U_{\tau^{\dag}}\le R$, we have $X_{\tau^{\dag}}\le R-1\le M-1$. Thus $\tau^{\dag}=\tau^{\rm o}$;
\item
if $X_0\ge M$, then $X_0,X_1,\dots,X_{\tau^--1}\ge0$ and $X_{\tau^-}<0$.
Since $X_{\tau^--1}=X_{\tau^-}-U_{\tau^-}$, $X_{\tau^-}\le -1$ and
$U_{\tau^-}\ge -L$, we have $X_{\tau^--1}\le L-1\le M-1$. Thus $\tau^{\rm o}\le\tau^--1$.
\end{itemize}

Now, we can observe the following connections between the foregoing times.
Invariance by translation implies that upshooting level $M$
(respectively downshooting level $-1$) when starting at a level $i\le M-1$
(respectively $i\ge 0$) is equivalent to upshooting level $0$
(respectively downshooting level $M-1$) when starting at a level $i-M$
(respectively $i+M$). In symbols, we have that
$$
\begin{array}{ll}
\mathbb{P}_i\{\tau^-=m,X_{\tau^-}=j\}=\mathbb{P}_{i+M}\{\tau^{\rm o}=m,X_{\tau^{\rm o}}=j+M\}
&\text{if $i\ge 0$,}
\\[1ex]
\mathbb{P}_i\{\tau^+=m,X_{\tau^+}=j\}=\mathbb{P}_{i-M}\{\tau^{\rm o}=m,X_{\tau^{\rm o}}=j-M\}
&\text{if $i\le M-1$,}
\end{array}
$$
and we deduce the identities below:
\begin{equation}\label{link-H-H+part}
H_{ij}^+(x)=H_{i-M\,j-M}^{\rm o}(x)\;\text{for any $i\le M-1$},\quad
H_{ij}^-(x)=H_{i+M\,j+M}^{\rm o}(x)\;\text{for any $i\ge 0$}.
\end{equation}
Additionally, by~(\ref{link-H-H+}),
\begin{equation}\label{link-H-Hdag}
H_{ij}^{\dag}(x)=\begin{cases}
H_{ij}^{\rm o}(x) & \text{if $i\le -1$,}
\\
p_{ij}\,x & \text{if $i\ge M$,}
\\[1ex]
\displaystyle x\Bigg(p_{ij}+\sum_{k=-M}^{-1} p_{ik} H_{kj}^{\rm o}(x)\Bigg)
& \text{if $0\le i\le M-1$.}
\end{cases}
\end{equation}

\subsection{Generating function of $X$}

Recall that $G_{ij}(x)=\sum_{m=0}^{\infty}\mathbb{P}_i\{X_m=j\}\,x^m.$
In the framework of random walk, we have the identity $G_{ij}=G_{0j-i}$
for any integers $i,j$;
then it is convenient to introduce the notation $\Gamma_{j}=G_{0j}$
so that $G_{ij}=\Gamma_{j-i}$.
We have the result below.
%
\begin{prop}\label{prop-Gij}
The generating function $\Gamma_{j-i}$ admits the following representation: for $x\in(-1,1)$,
\begin{equation}\label{*}
\Gamma_{j-i}(x)=\begin{cases}
\displaystyle\sum_{\ell\in\mathcal{L}^-}\frac{z_{\ell }(x)^{i-j+L-1}}{P'_x(z_{\ell}(x))}
&\text{if } i>j,
\\
\displaystyle -\sum_{\ell\in\mathcal{L}^+}\frac{z_{\ell }(x)^{i-j+L-1}}{P'_x(z_{\ell }(x))}
&\text{if } i\le j,
\end{cases}
\end{equation}
where $z_{\ell}(x)$, $\ell\in\{1,\dots,L+R\}$, are the roots of
the polynomial $P_x:z\mapsto z^L-x\sum_{j=0}^{L+R} \pi_{j-L} z^j$ and
$$
\mathcal{L}^+=\{\ell\in\{1,\dots,L+R\}: |z_{\ell}(x)|>1\},\quad
\mathcal{L}^-=\{\ell\in\{1,\dots,L+R\}: |z_{\ell}(x)|<1\}.
$$
\end{prop}
%
\Dim
We introduce the double generating function of the
numbers $\mathbb{P}_0\{X_m=j\}$, $m\in\mathbb{N},j\in\mathbb{Z}$:
$$
G(x,y)=\sum_{j\in\mathbb{Z}} \Gamma_j(x)\,y^j
=\sum_{m=0}^{\infty} \!\Bigg(\sum_{j\in\mathbb{Z}}\mathbb{P}_0\{X_m=j\}\,y^j\Bigg)x^m
=\sum_{m=0}^{\infty}\mathbb{E}_0\!\left(y^{X_m}\right)x^m
=\sum_{m=0}^{\infty}\left[x\,\mathbb{E}\!\left(y^{U_1}\right)\right]^m
$$
which can be simplified, for real numbers $x,y$ such that
$\left|x\,\mathbb{E}\!\left(y^{U_1}\right)\right|<1$, into
$$
G(x,y)=\frac{1}{1-x\,\mathbb{E}\!\left(y^{U_1}\right)}
=\frac{y^L}{y^L-x\sum_{j=0}^{L+R} \pi_{j-L}y^j}.
$$
Let us expand the rational fraction $y\mapsto G(x,y)$ into partial fractions.
For this, we introduce the polynomial $P_x(z)=z^L-x\sum_{j=0}^{L+R} \pi_{j-L} z^j$
and assume that $|x|<1$. We claim that the roots of $P_x$ have a modulus
different from $1$. Else, if $\mathrm{e}^{\mathrm{i}\theta}$ was a root of $P_x$,
we would have $\mathbb{E}\!\left(\mathrm{e}^{\mathrm{i}\theta U_1}\right)=1/x$.
This equality would simply entail that $|x|=1/\left|\mathbb{E}\!
\left(\mathrm{e}^{\mathrm{i}\theta U_1}\right)\right|\ge 1$ which contradicts our
assumption on $x$. Denote by $z_{\ell}(x)$, $\ell\in\{1,\dots,L+R\}$, the roots of $P_x$.
By the foregoing discussion, we can separate the roots having modulus greater
that 1 from those having modulus less that 1; they define the two sets
$\mathcal{L}^+$ and $\mathcal{L}^-$. 
With these settings at hand, we can write out the expansion of $G(x,y)$ as
$$
G(x,y)=\sum_{\ell=1}^{L+R} \frac{z_{\ell}(x)^L}{P'_x(z_{\ell}(x))}
\,\frac{1}{y-z_{\ell}(x)}.
$$
Next, by applying Taylor and Laurent series, we get
$$
\frac{1}{y-z_{\ell}(x)}
=\begin{cases}
\displaystyle -\sum_{j=0}^{+\infty}\frac{y^j}{z_{\ell}(x)^{j+1}}
&\text{if $\ell\in\mathcal{L}^+$ and $|y|<\min_{\ell\in\mathcal{L}^+}|z_{\ell}(x)|$,}
\\
\displaystyle\sum_{j=-\infty}^{-1}\frac{y^j}{z_{\ell}(x)^{j+1}}
&\text{if $\ell\in\mathcal{L}^-$ and $|y|>\max_{\ell\in\mathcal{L}^-}|z_{\ell}(x)|$,}
\end{cases}
$$
and, since $G_{ij}(x)=\Gamma_{j-i}(x)$, we extract by identification~(\ref{*}).
$\Box$

\subsection{Generating function of $T_n$ and $\tilde{T}_n$}\label{subsubsection-generating-Tn-LR}

Our aim is to apply Theorems~\ref{th-Kixy} and~\ref{th-Ktildeixy} to $(L,R)$-random walk.
First, we make the following observations which will be useful:
\begin{itemize}
\item
if $X_0\ge 0$, then $X_{\tau^-}\in\{-M,\dots,-1\}$ and if $X_0\le M-1$,
then $X_{\tau^{\dag}}\in\{0,\dots,2M-1\}$ and $X_{\tau^+}\in\{M,\dots,2M-1\}$;
\item
if $\tau^{\dag} >1$, then $X_{\tau^{\dag}}\le M-1$, or, equivalently,
if $X_{\tau^{\dag}}\ge M$, then $\tau^{\dag}=1$. Thus
$H_{ij}^{\dag}(y)=p_{ij} \,y$ if $j\ge M$.
\end{itemize}
Now, we rephrase Theorem~\ref{th-Kixy} in the present context.
%
\begin{teo}\label{th-Kixy-LR}
The generating functions $K_i$, $i\in\{0,\dots,M-1\}$, satisfy the system
\begin{align}
K_i(x,y)=\;
& x\sum_{j=0}^{M-1}
\!\left(p_{ij} +\sum_{k=M}^{2M-1} p_{ik} H_{kj}^{\rm o}(x)
+\sum_{k=-M}^{-1} p_{ik} H_{kj}^{\rm o}(y)\right)\! K_j(x,y)
+K_i(x,0)+K_i(0,y)-\mathbbm{1}_F(i)
\nonumber\\
& -x\sum_{j=0}^{M-1} \left(p_{ij} +\sum_{k=M}^{2M-1} p_{ik} H_{kj}^{\rm o}(x)\right)\! K_j(x,0),
\quad 0\le i\le M-1,
\label{eqKxy-LR}
\end{align}
where
\begin{align}
K_i(x,0)=\;
& G_i(x)-\sum_{j=0}^{M-1} H_{i+M\,j}^{\rm o}(x) G_{j-M}(x),
\nonumber\\[-2ex]
\label{eqKx0-K0y}\\[-2ex]
K_i(0,y)=\;
& G_i(y)-y\sum_{j=0}^{2M-1} p_{ij} G_j(y)
-\,y\sum_{j=0}^{M-1}\sum_{k=-M}^{-1} p_{ik} H_{kj}^{\rm o}(y) G_j(y),
\nonumber
\end{align}
and where the functions $H_{ij}^{\rm o}$ solve the systems
\begin{align}
\sum_{k=0}^{M-1} H_{ik}^{\rm o}(x)G_{kj}(x)=\;
& G_{ij}(x),\quad M\le i\le 2M-1,\; 0\le j\le M-1,
\nonumber
\\[-2ex]
\label{+1}
\\[-2ex]
\sum_{k=0}^{M-1} H_{ik}^{\rm o}(y)G_{kj}(y)=\;
& G_{ij}(y),\quad -M\le i\le -1,\; 0\le j\le M-1.
\nonumber
\end{align}
\end{teo}
%
\Dim
By the observations mentioned at the beginning of this subsubsection,
Equations~(\ref{eqK0xy}) take the form, for $0\le i\le M-1$,
\begin{align*}
K_i(x,0)=\; & G_i(x)-\sum_{j=-M}^{-1} H_{ij}^-(x) G_j(x),
\\
K_i(0,y)=\; & G_i(y)-\sum_{j=0}^{M-1} H_{ij}^{\dag}(y) G_j(y)
-\,y\sum_{j=M}^{2M-1}p_{ij} G_{j}(y),
\end{align*}
which can be rewritten, due to (\ref{link-H-H+part}) and~(\ref{link-H-Hdag}), as
(\ref{eqKx0-K0y}).
These latter contain the terms
$H_{ij}^{\rm o}(x)$, $M\le i\le 2M-1$, $0\le j\le M-1$, and
$H_{ij}^{\rm o}(y)$, $-M\le i\le -1$, $0\le j\le M-1$.
By~(\ref{cancel}) and (\ref{eqGH1}) we see that they solve Systems~(\ref{+1}).

Next, in order to apply~(\ref{eqKxy}) to compute $K_i(x,y)$,
we have to evaluate the quantities
$H_{ij}^{\rm o\dag}(x)$ and $H_{ij}^{\rm o-}(y)$. Because of~(\ref{cancel}),
in view of~(\ref{Ho-dag}) and~(\ref{Ho+-}), we have that
\begin{align*}
H_{ij}^{\rm o\dag}(x)=\; & x\Bigg(p_{ij} +
\sum_{k\in\mathbb{Z}:\atop M\le k\le i+R} p_{ik} H_{kj}^{\rm o}(x)\Bigg)\!
=x\Bigg(p_{ij} +\sum_{k=M}^{2M-1} p_{ik} H_{kj}^{\rm o}(x)\Bigg)\!,
\quad 0\le i,j\le M-1,
\\
H_{ij}^{\rm o-}(y)=\; & y\!\sum_{k\in\mathbb{Z}:\atop i-L\le k\le -1} p_{ik} H_{kj}^{\rm o}(y)
=y\sum_{k=-M}^{-1} p_{ik} H_{kj}^{\rm o}(y), \quad 0\le i,j\le M-1.
\end{align*}
Finally, putting these last equalities into~(\ref{eqKxy}) yields~(\ref{eqKxy-LR}).
$\Box$\\

Now, let us rephrase Theorem~\ref{th-Ktildeixy} in the present context.
Set $\varpi_i=\mathbb{P}_i\{X_1\le -1\text{ or } X_1\ge M\}$.
%
\begin{teo}\label{th-tilde-Kixy-LR}
The generating functions $\tilde{K}_i$, $i\in\{0,\dots,M-1\}$, satisfy the system
\begin{align}
\tilde{K}_i(x,y)=\;
& \sum_{j=0}^{M-1}\!\left(
x\sum_{k=M}^{2M-1}\sum_{\ell=0}^{M-1} p_{ik} H_{k\ell}^{\rm o}(x) H_{\ell j}^{\pm}(x)
+y\sum_{k=-M}^{-1}\sum_{\ell=0}^{M-1} p_{ik} H_{k\ell}^{\rm o}(y) H_{\ell j}^{\pm}(y)
\right)\! \tilde{K}_j(x,y)
\nonumber\\
& -x\sum_{j=0}^{M-1}\sum_{k=M}^{2M-1}\sum_{\ell=0}^{M-1}
p_{ik} H_{k\ell}^{\rm o}(x) H_{\ell j}^{\pm}(x) \tilde{K}_j(x,0)
-y\sum_{j=0}^{M-1}\sum_{k=-M}^{-1}\sum_{\ell=0}^{M-1}
p_{ik} H_{k\ell}^{\rm o}(y) H_{\ell j}^{\pm}(y) \tilde{K}_j(0,y),
\nonumber\\
& +\tilde{K}_i(x,0)+\tilde{K}_i(0,y)-\mathbbm{1}_F(i)\varpi_i,
\quad 0\le i\le M-1,
\label{eq-tilde-Kxy-LR}
\end{align}
where
\begin{align}
\tilde{K}_i(x,0)=\; &
x\!\sum_{j=M}^{2M-1}\! p_{ij} G_j(x)
-x\!\sum_{j=0}^{M-1}\sum_{k=2M}^{3M-1}\!  p_{i\,k-M} H_{kj}^{\rm o}(x) G_{j-M}(x)
+\mathbbm{1}_F(i)\varpi_i,
\nonumber
\\[-2ex]
\label{eq-tilde-Kx0-K0y}
\\[-2ex]
\tilde{K}_i(0,y)=\; & y\!\sum_{j=-M}^{-1}\! p_{ij} G_j(y)
-y\!\sum_{j=0}^{M-1}\!\sum_{k=-2M}^{-M-1}\!\!\!  p_{i\,k+M} H_{kj}^{\rm o}(y) G_{j+M}(y)
+\mathbbm{1}_F(i)\varpi_i,
\nonumber
\end{align}
and where the functions $H_{ij}^{\rm o}$ and $H_{ij}^{\pm}$ solve the systems
\begin{align}
&\sum_{k=0}^{M-1} H_{ik}^{\rm o}(x)G_{kj}(x)
=G_{ij}(x),\quad M\le i\le 2M-1 \text{ (resp. $2M\le i\le 3M-1$)},\; 0\le j\le M-1,
\label{+1bis}
\\
&\sum_{k=0}^{M-1} H_{ik}^{\rm o}(y)G_{kj}(y)
=G_{ij}(y),\quad -M\le i\le -1 \text{ (resp. $-2M\le i\le -M-1$)},\; 0\le j\le M-1,
\label{+2bis}
\\
&\sum_{k=0}^{M-1} (\delta_{ik}-p_{ik}x) H_{kj}^{\pm}(x)
=\delta_{ij}\varpi_i,
\quad 0\le i,j\le M-1.
\label{+3bis}
\end{align}
\end{teo}
%
\Dim
By the observations mentioned at the beginning of this subsubsection
and by applying formula~(\ref{G+-}), Equations~(\ref{tilde-K}) take the form, for $0\le i\le M-1$,
\begin{align}
\tilde{K}_i(x,0)=\; & G_i^+(x)-\sum_{j=-M}^{-1} H_{ij}^{-+}(x) G_j(x)
+\mathbbm{1}_F(i)\mathbb{P}_i\{X_1\le -1\},
\nonumber
\\[-2ex]
\label{eq-tilde-Kx0-K0y-inter}
\\[-2ex]
\tilde{K}_i(0,y)=\; & G_i^-(y)-\sum_{j=M}^{2M-1} H_{ij}^{+-}(y) G_j(y)
+\mathbbm{1}_F(i)\mathbb{P}_i\{X_1\ge M\}.
\nonumber
\end{align}
And, by~(\ref{cancel}), Equations~(\ref{H+-}) yield that
\begin{align}
H_{ij}^{-+}(x)=\;
&
x\sum_{k=M}^{2M-1} p_{ik} H_{kj}^-(x), \quad 0\le i\le M-1,\;-M\le j\le -1,
\nonumber
\\[-2ex]
\label{H+-LR}
\\[-2ex]
H_{ij}^{+-}(y)=\;
&
y\sum_{k=-M}^{-1} p_{ik} H_{kj}^+(y), \quad 0\le i\le M-1,\;M\le j\le 2M-1.
\nonumber
\end{align}
Actually, the terms
$H_{ij}^-(x)$, $M\le i\le 2M-1$, $-M\le j\le -1$ and
$H_{ij}^+(y)$, $-M\le i\le -1$, $M\le j\le 2M-1$
are directly related to the function $H_{ij}^{\rm o}$
according to~(\ref{link-H-H+part}). Hence, putting (\ref{link-H-H+part})
into (\ref{H+-LR}), and next the obtained equality
into Equations~(\ref{eq-tilde-Kx0-K0y-inter}),
these latter can be rewritten as (\ref{eq-tilde-Kx0-K0y}).

On the other hand,
by~(\ref{eqGH1}), the terms $H_{ij}^{\rm o}(x)$, $2M\le i\le 3M-1$, $0\le j\le M-1$, and
$H_{ij}^{\rm o}(y)$, $-2M\le i\le -M-1$, $0\le j\le M-1$ lying
in~(\ref{eq-tilde-Kx0-K0y}) solve Systems~(\ref{+1bis}) and (\ref{+2bis}).

Additionally, by rewriting Equations~(\ref{tilde-K*}) as
$$
\tilde{H}_{ij}^{\pm+}(x)
=x\sum_{k=M}^{2M-1}\sum_{\ell=0}^{M-1} p_{ik} H_{k\ell}^{\rm o}(x)H_{\ell j}^{\pm}(x),
\quad
\tilde{H}_{ij}^{\pm-}(y)
=y\sum_{k=-M}^{-1}\sum_{\ell=0}^{M-1} p_{ik} H_{k\ell}^{\rm o}(y)H_{\ell j}^{\pm}(y),
$$
and putting them into~(\ref{eqKxy-tilde}), we derive~(\ref{eq-tilde-Kxy-LR}).

Finally, thanks to~(\ref{eqHtildedag-}), we see that the terms
$H_{kj}^{\pm}$ solve System~(\ref{+3bis}).
$\Box$\\

Since all the systems displayed in Theorem~\ref{th-Kixy-LR}
and~\ref{th-tilde-Kixy-LR} are linear and include a finite number of equations,
it is natural to solve them by matrix calculus.
For this, we introduce the matrices below. We adopt the following intuitive
settings: the subscripts $\boldsymbol{+},\boldsymbol{\ddag}$ (double-plus),
$\boldsymbol{-},\boldsymbol{=}$ (double-minus) refer to the appearance
of the quantities $+M$, $+2M$, $-M$, $-2M$ in the generic term of the
corresponding matrices, the subscript $\boldsymbol{\dag}$ refers to the set of indices
$\{0,\dots,2M-1\}$ while the superscripts $\boldsymbol{+},\boldsymbol{-},
\boldsymbol{\pm}$ refer to the events $\{X_1\ge M\}$, $\{X_1\le -1\}$,
$\{X_1\le -1\text{ or }X_1\ge M\}$. So, with these conventions at hands, we set
$$
\mathbf{I}=\big(\delta_{ij}\big)_{0\le i,j\le M-1},\quad
\mathbf{I^{\boldsymbol{\pm}}}=\big(\delta_{ij}\varpi_i\big)_{0\le i,j\le M-1},\quad
\mathbf{1_F}=\big(\mathbbm{1}_F(i)\big)_{0\le i\le M-1},\quad
\mathbf{1_F^{\boldsymbol{\pm}}}=\big(\mathbbm{1}_F(i)\varpi_i\big)_{0\le i\le M-1},
$$
$$
\mathbf{P}=\big(\pi_{j-i}\big)_{0\le i,j\le M-1},\quad
\mathbf{\Gamma}(x)=\big(\Gamma_{j-i}(x)\big)_{0\le i,j\le M-1},
$$
$$
\mathbf{P_{\boldsymbol{-}}}=\big(\pi_{j-i-M}\big)_{0\le i,j\le M-1},\quad
\mathbf{P_{\boldsymbol{+}}}=\big(\pi_{j-i+M}\big)_{0\le i,j\le M-1},\quad
\mathbf{P_{\boldsymbol{\dag}}}=\big(\pi_{j-i}\big)_{0\le i\le M-1,0\le j\le 2M-1},
$$
$$
\mathbf{\Gamma_{\boldsymbol{-}}}(x)=\big(\Gamma_{j-i-M}(x)\big)_{0\le i,j\le M-1},\quad
\mathbf{\Gamma_{\boldsymbol{=}}}(x)=\big(\Gamma_{j-i-2M}(x)\big)_{0\le i,j\le M-1},
$$
$$
\mathbf{\Gamma_{\boldsymbol{+}}}(x)=\big(\Gamma_{j-i+M}(x)\big)_{0\le i,j\le M-1},\quad
\mathbf{\Gamma_{\boldsymbol{\ddag}}}(x)=\big(\Gamma_{j-i+2M}(x)\big)_{0\le i,j\le M-1},
$$
$$
\mathbf{G}(x)=\big(G_i(x)\big)_{0\le i\le M-1},\quad
\mathbf{G_{\boldsymbol{-}}}(x)=\big(G_{i-M}(x)\big)_{0\le i\le M-1},
$$
$$
\mathbf{G_{\boldsymbol{+}}}(x)=\big(G_{i+M}(x)\big)_{0\le i\le M-1},\quad
\mathbf{G_{\boldsymbol{\dag}}}(x)=\big(G_i(x)\big)_{0\le i\le 2M-1},
$$
$$
\mathbf{K}(x,y)=\big(K_i(x,y)\big)_{0\le i\le M-1},\quad
\mathbf{\tilde{K}}(x,y)=\big(\tilde{K}_i(x,y)\big)_{0\le i\le M-1}.
$$

%
\begin{ex}\label{example-M=2}
Consider the case where $M=2$, that is the case of the --at most-- four nearest neighbours
random walk (including the $(1,2)$-, $(2,1)$- and $(2,2)$-random walks). In this situation,
the walk is characterized by the non negative numbers $\pi_{-2},\pi_{-1},\pi_0,\pi_1,\pi_2$
such that $\pi_{-2}+\pi_{-1}+\pi_0+\pi_1+\pi_2=1$ and we have
$\pi_i=0$ for $i\in\mathbb{Z}\backslash\{-2,-1,0,1,2\}$.
Recall the notation $\Gamma_j=G_{0j}$ for any integer $j$.
Below, we rewrite the previous matrices:
$$
\mathbf{I}=\begin{pmatrix}1 & 0\\ 0 & 1\end{pmatrix}\!,\quad
\mathbf{I^{\boldsymbol{\pm}}}=\begin{pmatrix} 1-\pi_0-\pi_1 & 0\\
0 & 1-\pi_0-\pi_{-1}\end{pmatrix}\!,\quad
\mathbf{1_F}=\begin{pmatrix} \mathbbm{1}_F(0) \\ \mathbbm{1}_F(1) \end{pmatrix}\!,\quad
\mathbf{1_F^{\boldsymbol{\pm}}}=\begin{pmatrix} \mathbbm{1}_F(0)(1-\pi_0-\pi_1) \\
\mathbbm{1}_F(1)(1-\pi_0-\pi_{-1})\end{pmatrix}\!,
$$
$$
\mathbf{P}=\begin{pmatrix} \pi_0 & \pi_1 \\ \pi_{-1} & \pi_0\end{pmatrix}\!,\quad
\mathbf{P_{\boldsymbol{+}}}=\begin{pmatrix} \pi_2 & 0 \\ \pi_1 & \pi_2\end{pmatrix}\!,\quad
\mathbf{P_{\boldsymbol{-}}}=\begin{pmatrix} \pi_{-2} & \pi_{-1} \\ 0 & \pi_{-2}\end{pmatrix}\!,\quad
\mathbf{P_{\boldsymbol{\dag}}}=\begin{pmatrix} \pi_0 & \pi_1 & \pi_2 & 0 \\
\pi_{-1} & \pi_0 & \pi_1 & \pi_2\end{pmatrix}\!,
$$
$$
\mathbf{\Gamma}(x)
=\begin{pmatrix} \Gamma_0(x) & \Gamma_1(x)\\ \Gamma_{-1}(x) & \Gamma_0(x)\end{pmatrix}\!,\quad
\mathbf{\Gamma_{\boldsymbol{+}}}(x)
=\begin{pmatrix} \Gamma_2(x) & \Gamma_3(x)\\ \Gamma_1(x) & \Gamma_2(x)\end{pmatrix}\!,\quad
\mathbf{\Gamma_{\boldsymbol{-}}}(x)
=\begin{pmatrix} \Gamma_{-2}(x) & \Gamma_{-1}(x)\\ \Gamma_{-3}(x) & \Gamma_{-2}(x)\end{pmatrix}\!,
$$
$$
\mathbf{\Gamma_{\boldsymbol{=}}}(x)
=\begin{pmatrix} \Gamma_{-4}(x) & \Gamma_{-3}(x)\\ \Gamma_{-5}(x) & \Gamma_{-4}(x)\end{pmatrix}\!,\quad
\mathbf{\Gamma_{\boldsymbol{\ddag}}}(x)
=\begin{pmatrix} \Gamma_4(x) & \Gamma_5(x)\\ \Gamma_3(x) & \Gamma_4(x)\end{pmatrix}\!,
$$
$$
\mathbf{G}(x)=\begin{pmatrix} G_0(x) \\ G_1(x)\end{pmatrix}\!,\quad
\mathbf{G_{\boldsymbol{+}}}(x)=\begin{pmatrix} G_2(x) \\ G_3(x)\end{pmatrix}\!,\quad
\mathbf{G_{\boldsymbol{-}}}(x)=\begin{pmatrix} G_{-2}(x) \\ G_{-1}(x)\end{pmatrix}\!,\quad
\mathbf{G_{\boldsymbol{\dag}}}(x)=\begin{pmatrix} G_0(x) \\ G_1(x) \\ G_2(x) \\ G_3(x)\end{pmatrix}\!,
$$
$$
\mathbf{K}(x,y)=\begin{pmatrix} K_0(x,y) \\ K_1(x,y)\end{pmatrix}\!,\quad
\mathbf{\tilde{K}}(x,y)=\begin{pmatrix} \tilde{K}_0(x,y) \\ \tilde{K}_1(x,y)\end{pmatrix}\!.
$$
\end{ex}
%
\begin{teo}\label{th-K-Ktilde-LR}
The generating matrix $\mathbf{K}$ admits the representation
$\mathbf{K}(x,y)=\mathbf{D}(x,y)^{-1}\mathbf{N}(x,y)$ where
\begin{align}
\mathbf{D}(x,y)=\;
& \mathbf{I} -x\!\left(\mathbf{P}+\mathbf{P_{\boldsymbol{+}}}
\mathbf{\Gamma_{\boldsymbol{-}}}(x)\mathbf{\Gamma}(x)^{-1}
+\mathbf{P_{\boldsymbol{-}}}\mathbf{\Gamma_{\boldsymbol{+}}}(y)\mathbf{\Gamma}(y)^{-1}\right)\!,
\nonumber\\[-1ex]
\label{matrixK-LR-complete}
\\[-1ex]
\mathbf{N}(x,y)=\;
& \left(\mathbf{I} -x\big(\mathbf{P}+\mathbf{P_{\boldsymbol{+}}}
\mathbf{\Gamma_{\boldsymbol{-}}}(x)\mathbf{\Gamma}(x)^{-1}\big)
\right)\!\left(\mathbf{G}(x)-\mathbf{\Gamma_{\boldsymbol{-}}}(x)
\mathbf{\Gamma}(x)^{-1}\mathbf{G_{\boldsymbol{-}}}(x)\right)
\nonumber\\
& +\mathbf{G}(y)-y\,\mathbf{P_{\boldsymbol{\dag}}}\mathbf{G_{\boldsymbol{\dag}}}(y)
-y\,\mathbf{P_{\boldsymbol{-}}}\mathbf{\Gamma_{\boldsymbol{+}}}(y)\mathbf{\Gamma}(y)^{-1}
\mathbf{G}(y)-\mathbf{1_F}.
\nonumber
\end{align}
The generating matrix $\mathbf{\tilde{K}}$ admits the representation
$\mathbf{\tilde{K}}(x,y)=\mathbf{\tilde{D}}(x,y)^{-1}\mathbf{\tilde{N}}(x,y)$ where
\begin{align}
\mathbf{\tilde{D}}(x,y)=\;
& \mathbf{I}-x\,\mathbf{P_{\boldsymbol{+}}}\mathbf{\Gamma_{\boldsymbol{-}}}(x)
\mathbf{\Gamma}(x)^{-1}(\mathbf{I}-x\,\mathbf{P})^{-1} \mathbf{I^{\boldsymbol{\pm}}}
-y\,\mathbf{P_{\boldsymbol{-}}}\mathbf{\Gamma_{\boldsymbol{+}}}(y)\mathbf{\Gamma}(y)^{-1}
(\mathbf{I}-y\,\mathbf{P})^{-1} \mathbf{I^{\boldsymbol{\pm}}}
\vphantom{\mathbf{1_F^{\boldsymbol{\pm}}}},
\nonumber\\[-1ex]
\label{matrixK-tilde-LR-complete}
\\[-1ex]
\mathbf{\tilde{N}}(x,y)=\;
& \left(\mathbf{I}-x\,\mathbf{P_{\boldsymbol{+}}}\mathbf{\Gamma_{\boldsymbol{-}}}(x)
\mathbf{\Gamma}(x)^{-1}(\mathbf{I}-x\,\mathbf{P})^{-1} \mathbf{I^{\boldsymbol{\pm}}}\right)
\!\times\!\left(\mathbf{1_F^{\boldsymbol{\pm}}}
+x\,\mathbf{P_{\boldsymbol{+}}} \big(\mathbf{G_{\boldsymbol{+}}}(x)
-\mathbf{\Gamma_{\boldsymbol{=}}}(x) \mathbf{\Gamma}(x)^{-1}
\mathbf{G_{\boldsymbol{-}}}(x)\big)\right)
\nonumber\\
& +\left(\mathbf{I}-y\,\mathbf{P_{\boldsymbol{-}}}\mathbf{\Gamma_{\boldsymbol{+}}}(y)
\mathbf{\Gamma}(y)^{-1}(\mathbf{I}-y\,\mathbf{P})^{-1} \mathbf{I^{\boldsymbol{\pm}}}\right)
\!\times\!\left(\mathbf{1_F^{\boldsymbol{\pm}}}
+y\,\mathbf{P_{\boldsymbol{-}}} \big(\mathbf{G_{\boldsymbol{-}}}(y)
-\mathbf{\Gamma_{\boldsymbol{\ddag}}}(y) \mathbf{\Gamma}(y)^{-1}
\mathbf{G_{\boldsymbol{+}}}(y)\big)\right)-\mathbf{1_F^{\boldsymbol{\pm}}}.
\nonumber
\end{align}
\end{teo}
%

In Section~\ref{subsection-ordinary-RW}, we will apply this theorem to the
ordinary random walk (corresponding to the case $M=1$), and in
Subsubsection~\ref{subsection-22RW}, we will apply Theorem~\ref{th-K-Ktilde-LR}
to the case of the symmetric $(2,2)$-random walk.

\section{Ordinary random walk}\label{subsection-ordinary-RW}

In this section, we consider ordinary random walks on $\mathcal{E}=\mathbb{Z}$,
that is, those for which $L=R=1$ (and then $M=1$). They are characterized by the
three probabilities $\pi_{-1},\pi_0,\pi_1$ that we respectively relabel as $q,r,p$
for simplifying the settings:
$$
q=\mathbb{P}\{U_1=-1\},\quad r=\mathbb{P}\{U_1=0\},\quad p=\mathbb{P}\{U_1=1\}.
$$
We have of course $p+q+r=1$ and
$$
p_{ij}=\begin{cases}
q & \text{if $j=i-1$,}
\\
r & \text{if $j=i$,}
\\
p & \text{if $j=i+1$,}
\\
0 & \text{if $j\notin\{i-1,i,i+1\}$.}
\end{cases}
$$
We choose $E^{\rm o}=\{0\}$, $E^{\dag}=\{0,1,2,\dots\}$, $E^+=\{1,2,\dots\}$,
$E^-=\{\dots,-2,-1\}$. Without loss of generality, we focus on the
case where the starting point is located at $0$, i.e., in the current settings, we choose $i=0$.
Our aim is to apply Theorem~\ref{th-K-Ktilde-LR} to this example.

The roots of the polynomial $P_x(z)=-(px\,z^2-(1-rx)z+qx)$ are expressed as
$$
z(x)=\frac{1-rx-\sqrt{\Delta(x)}}{2px},\quad
\zeta(x)=\frac{1-rx+\sqrt{\Delta(x)}}{2px}
$$
where $\Delta(x)=(1-rx)^2-4pqx^2$. They are chosen such that
$|z(x)|<1<|\zeta(x)|$ for $x\in(0,1)$. We have $z(x)\zeta(x)=q/p$ and
$P'_x(z)=-2px\,z+(1-rx)$. In particular,
$P'_x(z(x))=-P'_x(\zeta(x))=\sqrt{\Delta(x)}$.

Formula~(\ref{*}) yields, for $x\in(0,1)$, that
$$
\Gamma_{j-i}(x)=\begin{cases}
\displaystyle\frac{z(x)^{i-j}}{\sqrt{\Delta(x)}} & \text{if $i\ge j$,}
\\[3ex]
\displaystyle\frac{\zeta(x)^{i-j}}{\sqrt{\Delta(x)}} & \text{if $i\le j$.}
\end{cases}
$$
In particular,
$$
\Gamma_0(x)=\frac{1}{\sqrt{\Delta(x)}},\quad
\Gamma_{-1}(x)=\displaystyle\frac{z(x)}{\sqrt{\Delta(x)}},\quad
\Gamma_1(x)=\frac{1}{\zeta(x)\sqrt{\Delta(x)}}
=\frac{pz(x)}{q\sqrt{\Delta(x)}}.
$$
$$
\Gamma_{-2}(x)=\displaystyle\frac{z(x)^2}{\sqrt{\Delta(x)}},\quad
\Gamma_2(x)=\frac{1}{\zeta(x)^2\sqrt{\Delta(x)}}
=\frac{p^2z(x)^2}{q^2\sqrt{\Delta(x)}}.
$$

Next, we rewrite the matrices of Subsubsection~\ref{subsubsection-generating-Tn-LR}
with the convention that for matrices with one entry,
we omit the parentheses and assimilate them to numbers:
$$
\mathbf{I}=1,\quad \mathbf{I^{\boldsymbol{\pm}}}
=\mathbf{1_F^{\boldsymbol{\pm}}}=\mathbbm{1}_F(0)(1-r),\quad
\mathbf{1_F}=\mathbbm{1}_F(0),\quad
\mathbf{P}=r,\quad \mathbf{P_{\boldsymbol{+}}}=p,\quad
\mathbf{P_{\boldsymbol{-}}}=q,\quad
\mathbf{P_{\boldsymbol{\dag}}}=\begin{pmatrix} r \;\; p\end{pmatrix}\!,
$$
$$
\mathbf{\Gamma}(x)=\Gamma_0(x) ,\quad
\mathbf{\Gamma_{\boldsymbol{-}}}(x)=\Gamma_{-1}(x),\quad
\mathbf{\Gamma_{\boldsymbol{+}}}(x)=\Gamma_1(x),\quad
\mathbf{\Gamma_{\boldsymbol{=}}}(x)=\Gamma_{-2}(x),\quad
\mathbf{\Gamma_{\boldsymbol{\ddag}}}(x)=\Gamma_2(x),
$$
$$
\mathbf{K}(x,y)=K_0(x,y),\quad \mathbf{\tilde{K}}(x,y)=\tilde{K}_0(x,y).
$$
For simplifying the forthcoming results, we will only consider the cases where
$F=\mathbb{Z}$ and $F=\{0\}$. In the case $F=\mathbb{Z}$, we have
$G_i(x)=1/(1-x)$ for any integer~$i$ and
$$
\mathbf{G}(x)=\mathbf{G_{\boldsymbol{+}}}(x)=
\mathbf{G_{\boldsymbol{-}}}(x)=\frac{1}{1-x},\quad
\mathbf{G_{\boldsymbol{\dag}}}(x)=\frac{1}{1-x}\begin{pmatrix} 1 \\ 1\end{pmatrix}\!,
$$
while in the case $F=\{0\}$, we have $G_i(x)=\Gamma_{-i}(x)$ for any integer~$i$ and
$$
\mathbf{G}(x)=\Gamma_0(x),\quad
\mathbf{G_{\boldsymbol{-}}}(x)=\Gamma_1(x),\quad
\mathbf{G_{\boldsymbol{+}}}(x)=\Gamma_{-1}(x),\quad
\mathbf{G_{\boldsymbol{\dag}}}(x)=\begin{pmatrix} \Gamma_0(x)\\ \Gamma_{-1}(x)\end{pmatrix}\!.
$$

Theorem~\ref{th-K-Ktilde-LR} provides the result below.
%
\begin{teo}\label{th-K-ordinary}
Set $\Delta(u)=(1-ru)^2-4pqu^2$. The generating function $K_0$
admits the following expression:
\begin{itemize}
\item
if $F=\mathbb{Z}$,
\begin{equation}\label{K0(x,y)}
K_0(x,y)=y\,\frac{(p-q)(x-y)+(1-y)\sqrt{\Delta(x)}+(1-x)\sqrt{\Delta(y)}}
{(1-x)(1-y)\big(y-x+y\sqrt{\Delta(x)}+x\sqrt{\Delta(y)}\,\big)};
\end{equation}
\item
if $F=\{0\}$,
\begin{equation}\label{K0(x,y)bis}
K_0(x,y)=\frac{2y}{y-x+y\sqrt{\Delta(x)}+x\sqrt{\Delta(y)}}.
\end{equation}
\end{itemize}
\end{teo}
%
\Dim
Observe that $\Gamma_{-1}(x)/\Gamma_0(x)=z(x)$ and $\Gamma_1(x)/\Gamma_0(x)
=1/\zeta(x)=pz(x)/q$.

\begin{itemize}
\item
In the case where $F=\mathbb{Z}$, in view of~(\ref{matrixK-LR-complete}), we see that
$K_0(x,y)$ can be written as $N(x,y)/D(x,y)$ where
\begin{align*}
(1-x)(1-y)N(x,y)=\;
& (1-y)\big(1-rx-px\,z(x)\big)\big(1-z(x)\big)
\\
& +(1-x)\big(1-(p+r)y-py\,z(y)\big)-(1-x)(1-y),
\\
D(x,y)=\;
& 1-rx-px\,z(x)-px\,z(y).
\end{align*}
Using the identity $px \,z(x)^2=(1-rx)z(x)-qx$, straightforward computations lead to
\begin{align*}
2(1-x)(1-y)N(x,y)=\;
& 2\big(1-(q+r)x-(p+r)y+rxy-px(1-y)z(x)-py(1-x)z(y)\big)
\\
=\; & (p-q)(x-y)+(1-y)\sqrt{\Delta(x)}+(1-x)\sqrt{\Delta(y)},
\\
2y\,D(x,y)=\;
& y-x+y\sqrt{\Delta(x)}+x\sqrt{\Delta(y)}.
\end{align*}
We immediately deduce~(\ref{K0(x,y)}).

\item
In the case where $F=\{0\}$, by~(\ref{matrixK-LR-complete}),
$K_0(x,y)$ can be written as $N(x,y)/D(x,y)$ where
$D(x,y)$ has exactly the same expression as previously and
$$
N(x,y)=\frac{1}{\sqrt{\Delta(x)}}\,\big(1-rx-px\,z(x)\big)\Big(1-\frac{p}{q}\,z(x)^2\Big)
+\frac{1}{\sqrt{\Delta(y)}}\,\big(1-ry-2py\,z(y)\big)-1.
$$
Thanks to $px \,z(x)^2=(1-rx)z(x)-qx$, elementary computations give
\begin{align*}
\big(1-rx-px\,z(x)\big)\Big(1-\frac{p}{q}\,z(x)^2\Big)=\;
& \frac{1}{2qx}\,\big(1-rx-pxz(x)\big)\big(2qx-(1-rx)z(x)\big)=\sqrt{\Delta(x)},
\\
1-ry-2py\,z(y)=\;
& \sqrt{\Delta(y)},
\end{align*}
so that $N(x,y)=1$ from which we extract~(\ref{K0(x,y)bis}).
\end{itemize}
$\Box$
%
\begin{prop}
For $F=\mathbb{Z}$, the following identity holds true:
\begin{equation}\label{K0(x,y)ter}
K_0(x,y)=K_0(x,0)K_0(0,y).
\end{equation}
\end{prop}
%
\Dim
Putting $x=0$ or $y\to 0$ into~(\ref{K0(x,y)})
and observing that $y-x+y\sqrt{\Delta(x)}+x\sqrt{\Delta(y)}\sim_{y\to 0}
y\big(1-rx+\sqrt{\Delta(x)}\,\big)$ entail that
$$
K_0(x,0)=\frac{1-(2q+r)x+\sqrt{\Delta(x)}}{(1-x)\big(1-rx+\sqrt{\Delta(x)}\,\big)}
\quad\mbox{and}\quad
K_0(0,y)=\frac{1-(2p+r)y+\sqrt{\Delta(y)}}{2(1-y)}.
$$
By elementary calculations, we see that
$$
\big((2p+r)x-1+\sqrt{\Delta(x)}\,\big)\big(1-rx+\sqrt{\Delta(x)}\,\big)
=2px\big(1-(2q+r)x+\sqrt{\Delta(x)}\,\big)
$$
and $K_0(x,0)$ can be simplified into
$$
K_0(x,0)=\frac{(2p+r)x-1+\sqrt{\Delta(x)}}{2px(1-x)}.
$$
Therefore, the product of the two terms $K_0(x,0)$ and $K_0(0,y)$ is given by
\begin{equation}\label{product}
K_0(x,0)K_0(0,y)=\frac{\big((2p+r)x-1+\sqrt{\Delta(x)}\,\big)
\big(1-(2p+r)y+\sqrt{\Delta(y)}\,\big)}{4px(1-x)(1-y)}.
\end{equation}
Straightforward but tedious computations that we do not report here
show that
\begin{align*}
\lqn{\big((2p+r)x-1+\sqrt{\Delta(x)}\,\big)\big(1-(2p+r)y+\sqrt{\Delta(y)}\,\big)
\big(y-x+y\sqrt{\Delta(x)}+x\sqrt{\Delta(y)}\,\big)}
=4pxy\big((p-q)(x-y)+(1-y)\sqrt{\Delta(x)}+(1-x)\sqrt{\Delta(y)}\,\big).
\end{align*}
Finally, comparing~(\ref{K0(x,y)}) and~(\ref{product}) ensures the validity of~(\ref{K0(x,y)ter}).
$\Box$\\

Decomposition~(\ref{K0(x,y)ter}) may be appropriate for inverting the
generating function $K_0$ in order to provide a closed form for the
probability distribution of $T_n$. We will not go further in this direction.

Concerning $\tilde{K}_0$, for $i=0$, we introduce the following
settings:
\begin{align*}
A(x,y)=\,
& (1-rx)(1-ry)\Big(2(1-r)(1-qx-py-rxy)
\\
&
-\frac{1}{px}(1-y)\big(r(1-rx)^2+(1-r)^2\,px\big)
-\frac{1}{qy}(1-x)\big(r(1-ry)^2+(1-r)^2\,qy\big)\Big),
\\
B(x,y)=\,
& \frac{1}{px}\,(1-y)(1-ry)\big(r(1-rx)^2+(1-r)^2\,px\big),
\\
C(x,y)=\,
& \frac{1}{qy}\,(1-x)(1-rx)\big(r(1-ry)^2+(1-r)^2\,qy\big).
\end{align*}
The quantity $A(x,y)$ admits the following expansion:
$$
A(x,y)=(1-rx)(1-ry)\Big(a_{11}xy+a_{10}x+a_{01}y+a_{1-1}\frac xy
+a_{-11}\frac yx-a_{-10}\frac{1}{x}-a_{0-1}\frac{1}{y}+a_{00}\Big),
$$
where the coefficients $a_{ij}$ (of $x^iy^j$) are expressed by means of $p,q,r$
as follows: by substituting $p+q=1-r$,
\begin{align*}
a_{11}&=\frac{2r^3}{p}+\frac{2r^3}{q}+2r^2-2r
=2r(1-r)\left(\frac{r^2}{pq}-2\right),
\\
a_{10}&=-\frac{r^3}{p}-\frac{2r^2}{q}+(1-r)^2-2q(1-r)
=-r^2\left(\frac{r}{p}+\frac{2}{q}\right)+(p-q)(1-r),
\\
a_{01}&=-\frac{r^3}{q}-\frac{2r^2}{p}+(1-r)^2-2p(1-r)
=-r^2\left(\frac{r}{q}+\frac{2}{p}\right)+(q-p)(1-r),
\\
a_{00}&=\frac{2r^2}{p}+\frac{2r^2}{q}-2r^2+2r
=2r(1-r)\left(\frac{r}{pq}+1\right),
\\
a_{-10}&=-\frac{r}{p},\quad
a_{0-1}=-\frac{r}{q},\quad
a_{1-1}=\frac{r}{q},\quad
a_{-11}=\frac{r}{p}.
\end{align*}
We also introduce
\begin{align*}
A'(x,y)&=(1-rx)(1-ry)\left(\big(1-2r-r^2\big)
+r\,\frac{(1-rx)^2}{2pqx^2}+r\,\frac{(1-ry)^2}{2pqy^2}\right),
\\
B'(x,y)&=-\frac{r(1-ry)}{pqx^2}\,\big(pq(1-r)x^2+(1-rx)^2\big),
\\
C'(x,y)&=-\frac{r(1-rx)}{pqy^2}\,\big(pq(1-r)y^2+(1-ry)^2\big).
\end{align*}

Theorem~\ref{th-K-Ktilde-LR} provides the result below.
%
\begin{teo}\label{th-Ktilde-ordinary}
The generating function $\tilde{K}_0$ admits the following expression:
\begin{itemize}
\item
if $F=\mathbb{Z}$,
\end{itemize}
\begin{equation}\label{K0(x,y)-tilde}
\tilde{K}_0(x,y)=\frac{1}{(1-x)(1-y)}\,
\frac{A(x,y)+B(x,y)\sqrt{\Delta(x)}+C(x,y)\sqrt{\Delta(y)}}
{2r(1-rx)(1-ry)+(1-r)\big((1-ry)\sqrt{\Delta(x)}+(1-rx)\sqrt{\Delta(y)}\,\big)};
\end{equation}
\begin{itemize}
\item
if $F=\{0\}$,
\end{itemize}
\begin{equation}\label{K0(x,y)-tildebis}
\tilde{K}_0(x,y)=\frac{A'(x,y)+B'(x,y)\sqrt{\Delta(x)}+C'(x,y)\sqrt{\Delta(y)}}
{2r(1-rx)(1-ry)+(1-r)\big((1-ry)\sqrt{\Delta(x)}+(1-rx)\sqrt{\Delta(y)}\,\big)}.
\end{equation}
\end{teo}
%
\Dim
\begin{itemize}
\item
In the case where $F=\mathbb{Z}$, in view of~(\ref{matrixK-tilde-LR-complete}),
we see that $\tilde{K}_0(x,y)$ can be written as $\tilde{N}(x,y)/\tilde{D}(x,y)$ where
\begin{align*}
\tilde{N}(x,y)=\;
& \Big(1-\frac{1-r}{1-rx}\,px\,z(x)\Big)\Big((1-r)+\frac{px}{1-x}\,\big(1-z(x)^2\big)\Big)
\\
& +\Big(1-\frac{1-r}{1-ry}\,py\,z(y)\Big)\Big((1-r)+\frac{qy}{1-y}\,
\big(1-\frac{p^2}{q^2}\,z(y)^2\big)\Big)-(1-r),
\\
\tilde{D}(x,y)=\;
& 1-\frac{1-r}{1-rx}\,px\,z(x)-\frac{1-r}{1-ry}\,py\,z(y).
\end{align*}
Using the identity $px \,z(x)^2=(1-rx)z(x)-qx$, we easily get that
\begin{align*}
\noalign{\hspace{2.2em} $(1-x)(1-y)(1-rx)(1-ry)\tilde{N}(x,y)$}
=\; & (1-y)(1-ry)\big(1-rx-(1-r)\,px \,z(x)\big)\big((1-r)-(1-rx)z(x)\big)
\\
& +(1-x)(1-rx)\big(1-ry-(1-r)\,py \,z(y)\big)
\Big((1-r)-\frac{p}{q}(1-ry)z(y)\Big)
\\
& -(1-r)(1-x)(1-y)(1-rx)(1-ry),
\\
\noalign{\hspace{2.2em} $(1-rx)(1-ry)\tilde{D}(x,y)$}
=\; & (1-rx)(1-ry)-(1-r)(1-ry)px\,z(x)-(1-r)(1-rx)py\,z(y).
\end{align*}
Straightforward computations lead to
\begin{align*}
2(1-x)(1-y)(1-rx)(1-ry)\tilde{N}(x,y)=\;
& A(x,y)+B(x,y)\sqrt{\Delta(x)}+C(x,y)\sqrt{\Delta(y)},
\\
2(1-rx)(1-ry)\tilde{D}(x,y)=\;
& 2r(1-rx)(1-ry)
+(1-r)\big((1-ry)\sqrt{\Delta(x)}+(1-rx)\sqrt{\Delta(y)}\,\big),
\end{align*}
from which we deduce~(\ref{K0(x,y)-tilde}).

\item
In the case where $F=\{0\}$, by~(\ref{matrixK-tilde-LR-complete}),
$\tilde{K}_0(x,y)$ can be written as $\tilde{N}(x,y)/\tilde{D}(x,y)$ where
$\tilde{D}(x,y)$ has exactly the same expression as previously and
\begin{align*}
\tilde{N}(x,y)=\;
& \Big(1-\frac{1-r}{1-rx}\,px\,z(x)\Big)\Big[(1-r)+\frac{px\,z(x)}{\sqrt{\Delta(x)}}
\Big(1-\frac{p}{q}\,z(x)^2\Big)\!\Big]
\\
& +\Big(1-\frac{1-r}{1-ry}\,py\,z(y)\Big)\Big[(1-r)+\frac{py\,z(y)}{\sqrt{\Delta(y)}}
\Big(1-\frac{p}{q}\,z(y)^2\Big)\!\Big]-(1-r)
\end{align*}
By replacing $z(x)$ by its expression, elementary computations give
\begin{align*}
\Big(1-\frac{1-r}{1-rx}\,px\,z(x)\Big)\Big[(1-r)+\frac{px\,z(x)}{\sqrt{\Delta(x)}}
\Big(1-\frac{p}{q}\,z(x)^2\Big)\!\Big]=\;
& \Big(1-\frac32\,r-\frac12 r^2+r\frac{(1-rx)^2}{2pqx^2}\Big)
\\
& -\frac{r}{2}\Big(\frac{1-rx}{pqx^2}+\frac{1-r}{1-rx}\Big)\sqrt{\Delta(x)},
\end{align*}
so that $N(x,y)=1$ from which we extract~(\ref{K0(x,y)-tildebis}).

\end{itemize}
$\Box$\\

In the case where $r=0$ (and $p+q=1$), that is, when the random walk does not stay
at its current location, the generating function $\tilde{K}_0$
can be simplified. We write its expression below.
%
\begin{cor}
In the case where $r=0$, we have $\Delta(x)=\sqrt{1-4pqx^2}$ and
the generating function $\tilde{K}_0$ is given by
\begin{itemize}
\item
if $F=\mathbb{Z}$,
$$
\tilde{K}_0(x,y)=\frac{(p-q)(x-y)+(1-y)\sqrt{\Delta(x)}+(1-x)\sqrt{\Delta(y)}}
{(1-x)(1-y)\big(\sqrt{\Delta(x)}+\sqrt{\Delta(y)}\,\big)};
$$
\item
if $F=\{0\}$,
$$
\tilde{K}_0(x,y)=\frac{1}{\sqrt{\Delta(x)}+\sqrt{\Delta(y)}}.
$$
\end{itemize}
\end{cor}
%
We retrieve the expressions obtained in~\cite{random-walk}; especially
in the case where $F=\mathbb{Z}$, the corresponding expression is rewritten in the form
$$
\tilde{K}_0(x,y)=\frac{1}{\sqrt{\Delta(x)}+\sqrt{\Delta(y)}}
\bigg(\frac{p-q+\sqrt{\Delta(x)}}{1-x}
+\frac{q-p+\sqrt{\Delta(y)}}{1-y}\bigg)
$$
which can be inverted in order to provide a closed form for the
probability distribution of $\tilde{T}_n$ (see \cite{random-walk}).

\section{Symmetric random walk}\label{subsection-symmetric}

In this part, we focus on the particular random walks satisfying $L=R=M$,
with steps lying in $\{-M,-M+1,\dots,M-1,M\}$, such that
$\pi_i=\pi_{-i}$ for all integer $i$. In this case, the random walk
$(X_m)_{m\in\mathbb{N}}$ is symmetric. We provide a representation
of the generating function of $\tau^{\rm o}$ which can be inserted into the matrices
introduced in Section~\ref{section-LR}.

\subsection{Generating function of $X$}

The polynomial $P_x$ takes the form
$$
P_x(z)=z^M\bigg[1-\pi_0x-x\sum_{j=1}^M \pi_j\!\left(z^j+\frac{1}{z^j}\right)\!\!\bigg]\!.
$$
To simplify the discussion, we will suppose throughout
Section~\ref{subsection-symmetric} that all the roots of $P_x$ are distinct.

We immediately see that its roots are inverse two by two so that the sets
$\mathcal{L}^+$ and $\mathcal{L}^-$ have the same cardinality $M$.
It is clear that there are $M$ roots of modulus less than one,
while their $M$ inverses have modulus greater than one.
So, we relabel the roots as $z_{\ell}(x)$ and
$z_{\ell+M}(x)=1/z_{\ell}(x)$,  $\ell\in\{1,\dots,M\}$
with the convention that $|z_{\ell}(x)|<1$. This yields
$\mathcal{L}^-=\{1,\dots,M\}$ and $\mathcal{L}^+=\{M+1,\dots,2M\}$.

Because of the equality $P_x(z)=z^{2M}P_x(1/z)$, we have
$P'_x(z)=2Mz^{2M-1}P_x(1/z)-z^{2M-2}P'_x(1/z)$ and this implies that
$$
\frac{z_{\ell}(x)^{M-1}}{P'_x(z_{\ell}(x))}
=-\frac{(1/z_{\ell}(x))^{M-1}}{P'_x(1/z_{\ell}(x))}.
$$
Therefore, expressions~(\ref{*}) of $\Gamma_{j-i}(x)$ can be simplified
and we get the generating function of $X$ which is displayed below.
%
\begin{prop}
The generating function of $X$ is characterized by the numbers
$G_{ij}(x)=\Gamma_{j-i}(x)$, $i,j\in\mathbb{Z}$, $x\in(0,1)$, with
\begin{equation}\label{**}
\Gamma_{j-i}(x)=\sum_{m=0}^{\infty}\mathbb{P}_i\{X_m=j\}\,x^m
=\sum_{\ell=1}^M\frac{z_{\ell}(x)^{M-1}}{P'_x(z_{\ell}(x))}\, z_{\ell}(x)^{|i-j|}.
\end{equation}
\end{prop}
%

Below, we give three examples of such random walks.
%
\begin{ex}
Let us consider the random walk with jump probabilities given by
$$
\begin{cases}
\displaystyle \pi_i=c\binom{2M}{i+M}\quad\text{for } i\in\{-M,-M+1,\dots,-1,1,\dots,M-1,M\},
\\[3ex]
\displaystyle \pi_0=1-c\left[4^M-\binom{2M}{M}\!\right]\!,
\end{cases}
$$
where $c$ is a positive constant such that
$c\le 1/\big[4^M-\binom{2M}{M}\big]$.
For $c=1/4^M$, we have $\pi_i=\binom{2M}{i+M}/4^M$
for any $i\in\{-M,-M+1,\dots,M-1,M\}$ and for $c=1/\big[4^M-\binom{2M}{M}\big]$,
we have $\pi_0=0$, that is, the walker never stays at its current position.

In this case,
$$
\mathbb{E}\!\left(y^{U_1}\right)\!=\frac{1}{y^M}
\left(c(y+1)^{2 M}+(1-c\,4^M)\,y^M\right),
\vspace{-.25\baselineskip}
$$
and
\vspace{-.25\baselineskip}
$$
G(x,y)=\frac{y^M}{P_x(y)}\quad\text{where}\quad
P_x(z)=\left(1-(1-c\,4^M)\,x\right)\! z^M-c\,x (z+1)^{2M}.
$$
Suppose that $x\in (0,1)$. The roots of $P_x$ are those of the $M$ quadratic equations
$$
(z+1)^2- e^{i\frac{2\pi}{M}r}
\sqrt{\frac{1-(1-c\,4^M)\,x}{c\,x}}\, z=0,\quad 0\le r\le M-1.
$$
We have that
\begin{align*}
P'_x(z_{\ell}(x))
&
=M\!\left(1-(1-c\,4^M)\,x\right)\! z_{\ell}(x)^{M-1}
-2Mc\,x (z_{\ell}(x)+1)^{2M-1}
=M\!\left(1-(1-c\,4^M)\,x\right)\!
\frac{1-z_{\ell}(x)}{1+z_{\ell}(x)}\,z_{\ell}(x)^{M-1}.
\end{align*}
Expression~(\ref{**}) takes the form
$$
\Gamma_{j-i}(x)=\frac{1}{M\left(1-(1-c\,4^M)\,x\right)}
\sum_{\ell=1}^M\frac{1+z_{\ell}(x)}{1-z_{\ell}(x)}\, z_{\ell}(x)^{|i-j|}.
$$
\end{ex}
%
\begin{ex}
Let us consider the random walk with jump probabilities given by
$$
\begin{cases}
\displaystyle \pi_i=c\rho^{|i|}\binom{M}{|i|}\quad\text{for }
i\in\{-M,-M+1,\dots,-1,1,\dots,M-1,M\},
\\[3ex]
\displaystyle \pi_0=1-2c\left((\rho+1)^M-1\right)\!,
\end{cases}
$$
where $c$ and $\rho$ are positive constants such that $c\le 1/\big(2(\rho+1)^M-1\big)$.
For $c=1/\big(2(\rho+1)^M-1\big)$, we have $\pi_0=0$, that is, the walker never
stays at its current position.

In this case,
$$
\mathbb{E}\!\left(y^{U_1}\right)\!=c\bigg((\rho y+1)^M+\left(\frac{\rho}{y}+1\right)^{\!M}
-2(\rho+1)^M\bigg)+1,
\vspace{-.25\baselineskip}
$$
and
\vspace{-.25\baselineskip}
$$
G(x,y)=\frac{y^M}{P_x(y)}\quad\text{where}\quad
P_x(z)=\left[1-\left(1-2c(\rho+1)^M\right)x\right]z^M
-cx\left[(\rho z^2+z)^M+(z+\rho)^M\right]\!.
$$
\end{ex}
%
\begin{ex}
Let us consider the random walk with jump probabilities given by
$$
\begin{cases}
\displaystyle \pi_i=c\quad\text{for }
i\in\{-M,-M+1,\dots,-1,1,\dots,M-1,M\},
\\[3ex]
\displaystyle \pi_0=1-2Mc,
\end{cases}
$$
where $c$ is a positive constant such that $c\le 1/(2M)$.
This is a random walk to the $2M$ nearest neighbours with identically
distributed jumps and a possible stay at the current position.
For $c=1/(2M)$, we have $\pi_0=0$, that is, the walker never stays at its current position.
For $c=1/(2M+1)$, each step put the walker to the $2M$ nearest neighbours
or let it at its current position with identical probability.

In this case,
$$
\mathbb{E}\!\left(y^{U_1}\right)\!=\big(1-(2M+1)c\big)+c\,\frac{1-y^{2M+1}}{y^M(1-y)},
\vspace{-.25\baselineskip}
$$
and
\vspace{-.25\baselineskip}
$$
G(x,y)=\frac{y^M}{P_x(y)}\quad\text{where}\quad
P_x(z)=\frac{1}{1-z}\left[\big(1-(1-(2M+1)c)x\big)\big(z^M-z^{M+1}\big)
-cx\big(1-z^{2M+1}\big)\right]\!.
$$
\end{ex}

\subsection{Symmetric $(2,2)$-random walk}\label{subsection-22RW}

In Example~\ref{example-M=2}, we have considered the case of non-symmetric
$(2,2)$-random walk. Now we have a look on the symmetric $(2,2)$-random walk
(corresponding to the case $M=2$)
which is characterized by the three non-negative numbers $\pi_0,\pi_1,\pi_2$
such that $\pi_0+2\pi_1+2\pi_2=1$:
$$
\pi_0=\mathbb{P}\{U_1=0\},\quad \pi_1=\mathbb{P}\{U_1=1\}=\mathbb{P}\{U_1=-1\},
\quad \pi_2=\mathbb{P}\{U_1=2\}=\mathbb{P}\{U_1=-2\}.
$$
We have that
$$
\mathbb{E}\!\left(y^{U_1}\right)\!=\pi_2\!\left(y^2+\frac{1}{y^2}\right)\!+
\pi_1\!\left(y+\frac{1}{y}\right)\!+\pi_0,
$$
and $G(x,y)=y^2/P_x(y)$ where
$$
P_x(z)=z^2-x\left(\pi_2z^4+\pi_1z^3+\pi_0z^2+\pi_1z+\pi_2\right)=
z^2\!\left[1-x\!\left(\pi_2\!\left(z^2+\frac{1}{z^2}\right)\!+
\pi_1\!\left(z+\frac{1}{z}\right)\!+\pi_0\right)\!\right]\!.
$$
Set $\delta(x)=(\pi_1+4\pi_2)^2+4\pi_2(1/x-1)$.
The roots of the polynomial $P_x$ can be explicitly calculated
by introducing the intermediate unknown $\zeta=z+1/z$ and solving
two quadratic equations. For $x\in(0,1)$, the roots are real and distinct,
those of absolute value less than $1$ are given by
\begin{align}
z_1(x)&=-\frac{1}{4\pi_2}\left(\pi_1-\sqrt{\delta(x)}
+\sqrt2\,\sqrt{\pi_1^2+4\pi_1\pi_2-2\pi_2+\frac{2\pi_2}{x}-\pi_1\sqrt{\delta(x)}}\right)\!,
\nonumber\\[-1ex]
\label{roots}
\\[-1ex]
z_2(x)&=-\frac{1}{4\pi_2}\left(\pi_1+\sqrt{\delta(x)}
+\sqrt2\,\sqrt{\pi_1^2+4\pi_1\pi_2-2\pi_2+\frac{2\pi_2}{x}+\pi_1\sqrt{\delta(x)}}\right)\!,
\nonumber
\end{align}
and the other ones are $z_3(x)=1/z_1(x)$ and $z_4(x)=1/z_2(x)$.
Moreover, rewriting $P_x$ as $P_x(z)=z^2Q_x(z+1/z)$ with
$Q_x(\zeta)=1-x\left(\pi_2\zeta^2+\pi_1\zeta+\pi_0-2\pi_2-1/x\right),$
we get that if $z$ is a root of $P_x$, then
$$
P'_x(z)=\left(z^2-1\right)Q'\!\left(z+\frac{1}{z}\right)\!
=x\left(1-z^2\right)\!\left(2\pi_2\!\left(z+\frac{1}{z}\right)\!+\pi_1\right)
$$
which can be simplified into
$$
P'_x(z_1(x))=x\left(1-z_1(x)^2\right)\sqrt{\delta(x)},\quad
P'_x(z_2(x))=-x\left(1-z_2(x)^2\right)\sqrt{\delta(x)}.
$$
Identity~(\ref{**}) yields the generating function of $X$ below.
%
\begin{prop}
The generating function of $X$ is characterized by the numbers
$G_{ij}(x)=\Gamma_{j-i}(x)$, $i,j\in\mathbb{Z}$, $x\in(0,1)$, with
$$
\Gamma_{j-i}(x)
=\frac{1}{x\sqrt{\delta(x)}}\left(\frac{z_1(x)^{|j-i|+1}}{1-z_1(x)^2}
-\frac{z_2(x)^{|j-i|+1}}{1-z_2(x)^2}\right)
$$
where $z_1(x)$ and $z_2(x)$ are displayed in~(\ref{roots}).
\end{prop}
%

Finally, we provide a representation of the generating functions of
sojourn times $T_n$ and $\tilde{T}_n$ in the case where $F=\mathbb{Z}$.
We apply Formulas~(\ref{matrixK-LR-complete}) and~(\ref{matrixK-tilde-LR-complete})
with $M=2$ and $F=\mathbb{Z}$.
%
\begin{teo}\label{th-gene-M=2}
The generating matrices of $T_n$ and $\tilde{T}_n$ admit the
respective representations
$$
\mathbf{K}(x,y)=\mathbf{D}(x,y)^{-1}\mathbf{N}(x,y)
\quad\text{and}\quad
\mathbf{\tilde{K}}(x,y)=\mathbf{\tilde{D}}(x,y)^{-1}\mathbf{\tilde{N}}(x,y)
$$
where $\mathbf{D}(x,y)^{-1},\mathbf{N}(x,y)$ are given by~(\ref{D})--(\ref{N}) further
and $\mathbf{\tilde{D}}(x,y)^{-1},\mathbf{\tilde{N}}(x,y)$ by~(\ref{Dtilde})--(\ref{Ntilde}).
\end{teo}
%
The explicit expressions of matrices $\mathbf{K}(x,y)$ and $\mathbf{\tilde{K}}(x,y)$
involving cumbersome computations, we postpone the proof of Theorem~\ref{th-gene-M=2}
to Subsection~\ref{proof-th-gene-M=2}.

\section{Proofs of Theorem~\ref{th-Kixy}, Theorem~\ref{th-Ktildeixy}
and Theorem~\ref{th-gene-M=2}}\label{proof-section}

In this section we give the proofs of Theorem~\ref{th-Kixy}, Theorem~\ref{th-Ktildeixy}
and all the  auxiliary results stated in Sections~\ref{section-settings}
and~\ref{section-gene-func} as well as the proof of Theorem~\ref{th-gene-M=2}.
Subsections~\ref{proof-Tn} and~\ref{proof-other-Tn} concern time $T_n$
while Subsections~\ref{proof-Am-Bm} and~\ref{proof-tilde-Tn}
concern time $\tilde{T}_n$; Subsection~\ref{proof-th-gene-M=2} concerns
the case of $(2,2)$-symmetric random walk.

\subsection{Proof of Theorem \ref{th-Kixy}}\label{proof-Tn}

We first observe that for $n\in\mathbb{N}^*$,
\begin{itemize}
\item
$T_n=0$ if and only if for all $m\in \{1,\dots, n\}$, $X_m\in E^-$, that is
if and only if $\tau^{\dag} > n$;
\item
$T_n=n$ if and only if for all $m\in\{1,\dots, n\}$, $X_m\in E^{\dag}$,
that is if and only if $\tau^->n$;
\item
$1\le T_n\le n-1$ if and only if there exists distinct integers
$\ell,\ell'\in\{1,\dots,n\}$ such that $X_{\ell}\in E^{\dag}$ and $X_{\ell'}\in E^-$. This is
equivalent to saying that $\tau^-\le n$ and $\tau^{\dag}\le n$.
\end{itemize}
In the last case, we have the three possibilities below:
\begin{itemize}
\item
if $X_0\in E^-$, then $\tau^{\rm o}=\tau^{\dag}\le n$ by Assumption~$(A_1)$;
\item
if $X_0\in E^+$, then $\tau^{\rm o}\le\tau^- -1\le n-1$ by Assumption~$(A_2)$;
\item
if $X_0\in E^{\rm o}$, the following possibilities occur:
    \begin{itemize}
    \item
    if $X_1\in E^{\rm o}$, then $\tau^{\rm o}=1$;
    \item
    if $X_1\in E^+$, then $\tau^{\rm o}\le\tau^- -1$. In this case,
    $\tau^{\rm o}$ is the return time to $E^{\rm o}$;
    \item
    if $X_1\in E^-$, then $\tau^{\rm o}=\tau^{\dag}$. In this case,
    $\tau^{\rm o}$ is the return time to $E^{\rm o}$.
    \end{itemize}
\end{itemize}
This discussion entails, for any $i\in\mathcal{E}$, that
\begin{align*}
\mathbb{P}_i\{T_n=0, X_n\in F\}=\;
& \mathbb{P}_i\{\tau^{\dag} >n, X_n\in F\}
=\mathbb{P}_i\{X_n\in F\}-\mathbb{P}_i\{\tau^{\dag}\le n, X_n\in F\}
\\
=\; & \mathbb{P}_i\{X_n\in F\}-\sum_{\ell=1}^n\mathbb{P}_i\{\tau^{\dag}=\ell, X_n\in F\}.
\end{align*}
From this, we deduce the following:
\begin{align*}
K_i(0,y)=\;
& \sum_{n=0}^{\infty}\mathbb{P}_i\{T_n=0, X_n\in F\}\,y^n
\\
=\;& \sum_{n=0}^{\infty}\mathbb{P}_i\{ X_n\in F\}\,y^n
-\sum_{n=1}^{\infty}\sum_{\ell=1}^n\mathbb{P}_i\{\tau^{\dag}=\ell, X_n\in F\}\,y^n
\\
=\; & \sum_{j\in F} \!\Bigg(\sum_{n=0}^{\infty}\mathbb{P}_i\{X_n=j\}\,y^n\Bigg)\!
-\sum_{\ell=1}^{\infty}\sum_{n=\ell}^{\infty}\sum_{j\in E^{\dag}}
\mathbb{P}_i\{\tau^{\dag}=\ell, X_{\tau^{\dag}}=j\}\,\mathbb{P}_j\{ X_{n-\ell}\in F\}\,y^n
\\
=\; & \sum_{j\in F} G_{ij}(y)-\sum_{j\in E^{\dag}}\Bigg(\sum_{\ell=1}^{\infty}
\mathbb{P}_i\{\tau^{\dag}=\ell, X_{\tau^{\dag}}=j\}\,y^{\ell}\Bigg)\!
\!\Bigg(\sum_{n=0}^{\infty}\mathbb{P}_j\{X_n\in F\}\,y^n\Bigg)
\\
=\; & G_i(y)-\sum_{j\in E^{\dag}} H_{ij}^{\dag}(y) G_j(y).
\end{align*}
This proves the second equality in~(\ref{eqK0xy}).
In a very similar way,
\begin{align*}
\mathbb{P}_i\{T_n=n, X_n\in F\}&=\mathbb{P}_i\{\tau^- > n, X_n\in F\}
=\mathbb{P}_i\{X_n\in F\}-\sum_{\ell=1}^n\mathbb{P}_i\{\tau^-=\ell, X_n\in F\}.
\end{align*}
This entails that
$$
K_i(x,0)=\sum_{n=0}^{\infty}\mathbb{P}_i\{T_n=n,X_n\in F\}\,x^n
=G_i(x)-\sum_{j\in E^-} H_{ij}^- (x) G_j(x)
$$
which proves the first equality in~(\ref{eqK0xy}).
Next, we compute $K_i(x,y)$:
\begin{align}
K_i(x,y)=\;& \sum_{n=0}^{\infty}\mathbb{P}_i\{T_n=n, X_n\in F\}\,x^n
+\sum_{n=0}^{\infty}\mathbb{P}_i\{T_n=0, X_n\in F\}\,y^n -\mathbb{P}_i\{X_0\in F\}
\nonumber\\
& +\!\!\sum_{m,n\in\mathbb{N}:\atop 1\le m\le n-1}
\mathbb{P}_i\{T_n=m, X_n\in F\}\,x^m y^{n-m}
\nonumber\\
=\;& K_i(x,0)+K_i(0,y)+\!\!\sum_{m,n\in\mathbb{N}:\atop 1\le m\le n-1}
\mathbb{P}_i\{\tau^{\rm o}\le n,T_n=m, X_n\in F\}\,x^m y^{n-m}-\mathbbm{1}_F(i).
\label{Ki(x,y)}
\end{align}
We observe that on the set $\{\tau^{\rm o}\le n\}$
$$
T_n=\sum_{m=1}^n\mathbbm{1}_{E^{\dag}}(X_m)
=T_{\tau^{\rm o}}+\sum_{m=\tau^{\rm o}+1}^n\mathbbm{1}_{E^{\dag}}(X_m)
=T_{\tau^{\rm o}}+\sum_{m=1}^{n-\tau^{\rm o}}\mathbbm{1}_{E^{\dag}}(X_{m+\tau^{\rm o}})
=T_{\tau^{\rm o}}+T^{(\tau^{\rm o})}_{n-\tau^{\rm o}}
$$
where $T^{(\tau^{\rm o})}_{n-\tau^{\rm o}}$ is the sojourn time in $E^{\dag}$ up to time $n-\tau^{\rm o}$
for the shifted chain $(X_{m+\tau^{\rm o}})_{m\in\mathbb{N}}$. Moreover,
\begin{itemize}
\item
if $X_1\in E^{\dag}$, we have $X_1,\dots,X_{\tau^{\rm o}}\in E^{\dag}$ and then
$T_{\tau^{\rm o}}=\tau^{\rm o}$ which yields
$T_{n-\tau^{\rm o}}^{(\tau^{\rm o})}=T_n-\tau^{\rm o}$ ;
\item
if $X_1\in E^-$, we have $X_1,\dots,X_{\tau^{\rm o}-1}\in E^-,X_{\tau^{\rm o}}\in E^{\rm o}$;
in this case $T_{\tau^{\rm o}}=1$ and then $T_{n-\tau^{\rm o}}^{(\tau^{\rm o})}=T_n-1$.
\end{itemize}
Consequently, for $m\in\{1,\dots,n-1\}$ and $i\in\mathcal{E}$,
\begin{align*}
\mathbb{P}_i\{\tau^{\rm o} \le n, T_n=m, X_n\in F\}=\;
& \mathbb{P}_i\{X_1 \in E^{\dag},\tau^{\rm o} \le n, T_n=m, X_n\in F\}
\\
& +\mathbb{P}_i\{X_1 \in E^-,\tau^{\rm o} \le n, T_n=m, X_n\in F\}
\\
=\;& \sum_{j \in E^{\rm o}} \sum_{\ell=1}^{n-1}
\mathbb{P}_i\{X_1\in E^{\dag},\tau^{\rm o}=\ell, X_{\tau^{\rm o}}=j\}\,
\mathbb{P}_j\{T_{n-\ell}=m-\ell, X_{n-\ell}\in F\}
\\
& +\sum_{j \in E^{\rm o}} \sum_{\ell=2}^{n}
\mathbb{P}_i\{X_1\in E^-,\tau^{\rm o}=\ell, X_{\tau^{\rm o}}=j\}\,
\mathbb{P}_j\{T_{n-\ell}=m-1, X_{n-\ell}\in F\}.
\end{align*}
The sum in~(\ref{Ki(x,y)}) can be evaluated as follows:
\begin{align}
\lqn{
\sum_{m,n\in\mathbb{N}:\atop 1\le m\le n-1}
\mathbb{P}_i\{\tau^{\rm o}\le n, T_n=m, X_n\in F\}\,x^m y^{n-m}}=\;
& \sum_{\ell=1}^{\infty} \sum_{j\in E^{\rm o}}
\mathbb{P}_i\{X_1\in E^{\dag},\tau^{\rm o}=\ell, X_{\tau^{\rm o}}=j\}
\sum_{m,n\in\mathbb{N}:\atop\ell\le m\le n-1}
\mathbb{P}_j\{T_{n-\ell}=m-\ell, X_{n-\ell}\in F\}\,x^{m} y^{n-m}
\nonumber\\
& + \sum_{\ell=2}^{\infty} \sum_{j\in E^{\rm o}}
\mathbb{P}_i\{X_1\in E^-,\tau^{\rm o}=\ell, X_{\tau^{\rm o}}=j\}
\sum_{m,n\in\mathbb{N}:\atop 1\le m\le n-\ell+1 }
\mathbb{P}_j\{ T_{n-\ell}=m-1, X_{n-\ell}\in F\}\,x^{m} y^{n-m}
\nonumber\\
=\; &\sum_{j\in E^{\rm o}}\sum_{\ell=1}^{\infty}
\mathbb{P}_i\{X_1\in E^{\dag},\tau^{\rm o}=\ell, X_{\tau^{\rm o}}=j\}\,x^{\ell}
\sum_{m,n\in\mathbb{N}:\atop m\le n-1 }
\mathbb{P}_j\{ T_{n}=m, X_{n}\in F\}\,x^{m} y^{n-m}
\nonumber\\
& +\sum_{j\in E^{\rm o}}\sum_{\ell=2}^{\infty}
\mathbb{P}_i\{X_1\in E^-,\tau^{\rm o}=\ell, X_{\tau^{\rm o}}=j\}\,x y^{\ell-1}
\sum_{m,n\in\mathbb{N}:\atop m\le n }
\mathbb{P}_j\{ T_{n}=m, X_{n}\in F\}\,x^{m} y^{n-m}
\nonumber\\
=\;& \sum_{j\in E^{\rm o}}\mathbb{E}_i \big(x^{\tau^{\rm o}}
\mathbbm{1}_{\{X_1\in E^{\dag},X_{\tau^{\rm o}}=j\}}\big) [K_j (x,y)-K_j(x,0)]
+\frac x y\sum_{j\in E^{\rm o}}\mathbb{E}_i \big(y^{\tau^{\rm o}}
\mathbbm{1}_{\{X_1\in E^-,X_{\tau^{\rm o}}=j\}}\big) K_j(x,y)
\nonumber\\
=\;& \sum_{j\in E^{\rm o}} \!\left(H_{ij}^{\rm o\dag} (x)+\frac x y\,H_{ij}^{\rm o-} (y)\right)\!
K_j(x,y)-\sum_{j\in E^{\rm o}} H_{ij}^{\rm o\dag} (x) K_j(x,0).
\label{sum-in-Ki(x,y)}
\end{align}
As a result, putting~(\ref{sum-in-Ki(x,y)}) into~(\ref{Ki(x,y)}), we get~(\ref{eqKxy}).
$\Box$

\subsection{Proof of Proposition \ref{propKixy}}\label{proof-other-Tn}

By putting~(\ref{eqK0xy}) into~(\ref{eqKxy}), we obtain, for all $i\in\cal{E}$, that
\begin{align}
K_i(x,y)=\;
& G_i(x)-\sum_{j\in E^-} H_{ij}^-(x) G_j(x)+G_i(y)
-\sum_{j\in E^{\dag}} H_{ij}^{\dag}(y) G_j(y)
-\sum_{j\in E^{\rm o}} H_{ij}^{\rm o\dag}(x) G_j(x)
-\mathbbm{1}_F(i)
\nonumber\\
& +\sum_{j\in E^{\rm o}}\!\left(H_{ij}^{\rm o\dag}(x)
+\frac x y\,H_{ij}^{\rm o-}(y)\right)\! K_j(x,y)
+\sum_{k\in E^{\rm o}} H_{ik}^{\rm o\dag}(x)\sum_{j\in E^-} H_{kj}^-(x) G_j(x).
\label{eqKxy-bis}
\end{align}
If $i\in E^+$, we have $H_{ij}^{\dag}(y)=p_{ij}\,y$,
$H_{ij}^{\rm o\dag}(x)=H_{ij}^{\rm o}(x)$, and $H_{ij}^{\rm o-}(y)=0$. Putting this into
(\ref{eqKxy-bis}), we get that
\begin{align}
K_i(x,y)=\;
& G_i(x)-\sum_{j\in E^{\rm o}} H_{ij}^{\rm o}(x) G_j(x)
+\sum_{j\in E^{\rm o}} H_{ij}^{\rm o}(x) K_j(x,y)
\nonumber\\
& +\sum_{j\in E^-} \!\Bigg(\sum_{k\in E^{\rm o}} H_{ik}^{\rm o}(x) H_{kj}^-(x)-H_{ij}^-(x)
\Bigg)G_j(x)+ \!\Bigg(G_i(y)-y\sum_{j\in E^{\dag}} p_{ij} G_j(y)-\mathbbm{1}_F(i)\Bigg)\!.
\label{eqKxy-ter}
\end{align}
We claim that both terms within brackets in the second line of~(\ref{eqKxy-ter})
vanish. Indeed, for $j\in E^-$, since $\tau^{\rm o}<\tau^-$, we have that
$$
H_{ij}^-(x)=\mathbb{E}_i \big(x^{\tau^-}\mathbbm{1}_{\{X_{\tau^-}=j\}}\big)
=\sum_{k\in E^{\rm o}}\mathbb{E}_i\big(x^{\tau^{\rm o}} \mathbbm{1}_{\{X_{\tau^{\rm o}}=k\}}\big)\,
\mathbb{E}_k \big(x^{\tau^-} \mathbbm{1}_{\{X_{\tau^-}=j\}}\big)
=\sum_{k\in E^{\rm o}} H_{ik}^{\rm o}(x)H_{kj}^-(x).
$$
Therefore, the first term within brackets in the second line of~(\ref{eqKxy-ter})
vanishes. Moreover, by~(\ref{eqG1}), we have that
$$
G_i(y)=\sum_{k\in F} G_{ik}(y)
=\sum_{k\in F} \!\Bigg(\delta_{ik} +y\sum_{j\in E^{\dag}} p_{ij} G_{jk}(y)\Bigg)\!
=y\sum_{j\in E^{\dag}} p_{ij} G_j(y)+\mathbbm{1}_F(i)
$$
which proves that the second term within brackets in~(\ref{eqKxy-ter})
vanishes too. Hence we have checked that~(\ref{uno}) holds true.

On the other hand, for $i\in E^-$, we have necessarily $\tau^{\rm o}=\tau^{\dag}$ and
then $H_{ij}^{\dag}(y)=H_{ij}^{\rm o}(y)\,\mathbbm{1}_{E^{\rm o}} (j)$.
Moreover, if  $i\in E^-$ and $X_1\in E^{\dag}$, then $\tau^{\rm o}=1$ and
$X_1\in E^{\rm o}$ and we have that
$$
H_{ij}^{\rm o\dag}(x)=p_{ij}\,x\,\mathbbm{1}_{E^{\rm o}}(j),
\quad H_{ij}^{\rm o-}(y)=\mathbb{E}_i\big(y^{\tau^{\rm o}}\mathbbm{1}_{\{X_{\tau^{\rm o}}=j\}}\big)
-\mathbb{E}_i\big(y^{\tau^{\rm o}} \mathbbm{1}_{\{X_1\in E^{\dag},X_{\tau^{\rm o}}=j\}}\big)
=H_{ij}^{\rm o}(y)- p_{ij}\,y\,\mathbbm{1}_{E^{\rm o}}(j).
$$
Then, putting these equalities into~(\ref{eqKxy-bis}), we get that
\begin{align}
K_{i}(x,y)=\;
& G_i(y)-\sum_{j\in E^{\rm o}} H_{ij}^{\rm o}(y) G_j(y)
+\frac x y\sum_{j\in E^{\rm o}}H_{ij}^{\rm o}(y) K_j(x,y)
\nonumber\\
& +\!\Bigg(G_i(x)-x\sum_{j\in E^{\rm o}} p_{ij} G_{j}(x)-\mathbbm{1}_F(i)\Bigg)\!
+\sum_{j\in E^-} \!\Bigg(x\sum_{k\in E^{\rm o}} p_{ik}  H_{kj}^-(x)-
H_{ij}^-(x)\Bigg)G_j(x).
\label{eqKxy-quater}
\end{align}
We claim that the sum of both terms within brackets in the second line
of~(\ref{eqKxy-quater}) vanish. This can be easily seen by using~(\ref{eqG1})
and Assumption~$(A_1)$  which yield that
$$
G_i(x)-x\sum_{j\in E^{\rm o}} p_{ij} G_j(x)-\mathbbm{1}_F(i)
=x\sum_{j\in E^-} p_{ij} G_j(x)
$$
and by observing that
$$
H_{ij}^-(x)=\mathbb{E}_i\big(x^{\tau^-}\mathbbm{1}_{\{X_1\in E^-,X_{\tau^-}=j\}}\big)
+\mathbb{E}_i\big(x^{\tau^-}\mathbbm{1}_{\{X_1\in E^{\rm o},X_{\tau^-}=j\}}\big)
=x\sum_{k\in E^{\rm o}} p_{ik} H_{kj}^-(x)+p_{ij}\,x\,\mathbbm{1}_{E^-}(j).
$$
Finally, Equality~(\ref{eqKxy-quater}) can be simplified to~(\ref{due}).
$\Box$

\subsection{Proof of Proposition \ref{lemma}}\label{proof-Am-Bm}

Set $A_{0,m}^{\rm o}=\{X_m\in E^{\rm o}\}$.
We notice that for any $m\in\mathbb{N}^*$ and any $\ell\in\{1,\dots,m\}$,
$A_{\ell,m}=A_{0,m}^{\rm o}\cap A_{\ell-1,m-1}$. Therefore, putting this into
the definition of $B_m$, we get, for any $m\in\mathbb{N}^*$, that
\begin{align*}
B_m=\;
&A_{0,m}\cup A_{1,m}\cup\dots\cup A_{m-1,m}
\\
=\;& A_{0,m}\cup \left(A_{0,m}^{\rm o}\cap A_{0,m-1}\right)\cup\dots
\cup \left(A_{0,m}^{\rm o}\cap A_{m-2,m-1}\right)
\\
=\;& A_{0,m}\cup \left[A_{0,m}^{\rm o}\cap(A_{0,m-1}\cup \dots\cup A_{m-2,m-1})\right]
\end{align*}
and we see, for any $m\in\mathbb{N}\backslash\{0,1\}$, that
\begin{equation}\label{Bm}
B_m=A_{0,m}\cup \left(A_{0,m}^{\rm o}\cap B_{m-1}\right)\!.
\end{equation}
In the same way, we have that
\begin{equation}\label{Bm'}
B'_m=A'_{0,m}\cup \left(A_{0,m}^{\rm o}\cap B'_{m-1}\right)\!.
\end{equation}
By observing that $\left(A_{0,m}^{\rm o}\right)^c\cap A_{0,m}^c=A'_{0,m}$
and using the elementary equality $A\cup B=A\cup(B\backslash A)$, we obtain
from~(\ref{Bm}) that
\begin{align*}
B_m^c=\;
& A_{0,m}^c\cap \left(\left(A_{0,m}^{\rm o}\right)^c\cup B_{m-1}^c\right)
\\
=\; & \left(\left(A_{0,m}^{\rm o}\right)^c\cap A_{0,m}^c\right)\cup
\left(A_{0,m}^c\cap B_{m-1}^c\right)
\\
=\; & A'_{0,m}\cup\left[\left(A_{0,m}^c\cap B_{m-1}^c\right)\backslash A'_{0,m}\right]\!.
\end{align*}
The term within brackets in the foregoing equality can be simplified as follows:
$$
\left(A_{0,m}^c\cap B_{m-1}^c\right)\backslash A'_{0,m}
=\left(A_{0,m}^c\backslash A'_{0,m}\right)\cap B_{m-1}^c
=A_{0,m}^{\rm o}\cap B_{m-1}^c.
$$
As a by-product, for any $m\in\mathbb{N}^*$,
\begin{equation}\label{Bmc}
B_m^c=A'_{0,m}\cup \left(A_{0,m}^{\rm o}\cap B_{m-1}^c\right)\!.
\end{equation}
Now, in view of~(\ref{Bm'}) and~(\ref{Bmc}), it is clear by recurrence
that $B_m^c=B'_m$ as claimed.
$\Box$

\subsection{Proof of Theorem~\ref{th-Ktildeixy}}\label{proof-tilde-Tn}

The observations below will be useful to compute the generating function
of $\tilde{T}_n$.
%
\begin{lem}\label{lemma-discussion}
Assume that $X_0\in E^{\rm o}$. The following hold true: for any $n\in\mathbb{N}^*$,
\begin{itemize}
\item
$\tilde{T}_n=0$ if and only if $\tau^+>n$;
\item
$\tilde{T}_n=n$ if and only if $X_1\in E^+$ and $\tau^->n$;
\item
if $1\le\tilde{T}_n\le n-1$, then $\tilde{\tau}^{\pm}\le n$ and
$\tilde{T}_n=(\tilde{\tau}^{\pm}-1)\mathbbm{1}_{E^+}+\tilde{T}'_n$
where $\tilde{T}'_n=\sum_{m=\tilde{\tau}^{\pm}}^n \delta_m$.
Additionally, $\tilde{T}'_n$ has the same distribution as
$\tilde{T}_{n-\tilde{\tau}^{\pm}+1}$ and is independent from $X_1,\dots,X_{\tilde{\tau}^{\pm}-1}$.
\end{itemize}
\end{lem}
%
\Dim
Fix a positive integer $n$.
\begin{itemize}
\item
Concerning the first point, we have that
    \begin{itemize}
    \item
    the set $\{X_1,\dots,X_m\in E^-\cup E^{\rm o}\}$ is included
    in $B'_m$ for all $m\in\{1,\dots,n\}$. Then, if $X_m\in E^-\cup E^{\rm o}$
    for all $m\in\{1,\dots,n\}$, we have $\delta_1=\dots=\delta_n=0$ and $\tilde{T}_n=0$;
    \item
    if there existed an $\ell\in\{1,\dots, n\}$ such that
    $X_{\ell}\in E^+$, we would have of course $\tilde{T}_n\ge 1$.
    \end{itemize}
This is equivalent to saying that $\tau^+>n$.

\item
Concerning the second point, we have that
    \begin{itemize}
    \item
    the set $\{X_1\in E^+,X_2,\dots,X_m\in E^+\cup E^{\rm o}\}$
    is included in $B_m$ for all $m\in\{1,\dots,n\}$. Then,
    if $X_1\in E^+$ and $X_2,\dots,X_n\in E^+\cup E^{\rm o}$,
    we have $\delta_1=\dots=\delta_n=1$ and $\tilde{T}_n=n$;
    \item
    if there existed an $\ell\in\{2,\dots, n\}$ such that
    $X_{\ell}\in E^-$, we would have of course $\tilde{T}_n\le n-1$.
    If we had $X_1\in E^-\cup E^{\rm o}$, we would have $\delta_1=0$ and
    $\tilde{T}_n\le n-1$ too.

    \end{itemize}
This is equivalent to saying that $X_1\in E^+$ and $\tau^->n$.

\item
Concerning the third point, we have that
    \begin{itemize}
    \item
    $1\le\tilde{T}_n\le n-1$ if and only if there exists distinct
    integers $\ell,\ell'\in\{1,\dots,n\}$ such that $\delta_{\ell}=1$ and $\delta_{\ell'}=0$.
    This implies that there exists necessarily an integer $m\in\{1,\dots,n\}$ such that
    $X_m\in E^{\rm o}$. Else, because of Assumptions~$(A_1)$ and $(A_2)$, we would have
    $X_m\in E^+$ for all $m\in\{1,\dots,n\}$
    or $X_m\in E^-$ for all $m\in\{1,\dots,n\}$;
    in both cases, we would have $\tilde{T}_n\in\{0,n\}$. Therefore $1\le\tau^{\rm o}\le n$.
    \item
    On the other hand, we must take care of the successive points lying
    in $E^{\rm o}$ after time $\tau^{\rm o}$. For describing them, time $\tilde{\tau}^{\pm}$
    will be useful.
    Obviously, we have $\tilde{\tau}^{\pm}>\tau^{\rm o}\ge 1$
    and between times $\tau^{\rm o}$ and $\tilde{\tau}^{\pm}-1$, the chain
    $(X_n)_{n\in\mathbb{N}}$ stays in $E^{\rm o}$.
    More precisely, we must take care of the points corresponding to
    the times between $\tau^{\rm o}$ and $(\tilde{\tau}^{\pm}-1)\wedge n$.
    In fact, these $(\tilde{\tau}^{\pm}-1)\wedge n-\tau^{\rm o}$ points are counted in $\tilde{T}_n$
    if $X_1\in E^+$, not counted if $X_1\in E^-\cup E^{\rm o}$.
    Then, if $\tilde{\tau}^{\pm}\le n$, we have $\tilde{T}_n=\tilde{T}_{\tilde{\tau}^{\pm}-1}
    +\tilde{T}'_n$ with
    $$
    \tilde{T}_{\tilde{\tau}^{\pm}-1}=\begin{cases}
    \tilde{\tau}^{\pm}-1 &\text{if $X_1\in E^+$,}
    \\
    0 & \text{if $X_1\in E^-\cup E^{\rm o}$,}
    \end{cases}
    \quad\text{and}\quad
    \tilde{T}'_n=\sum_{m=\tilde{\tau}^{\pm}}^n \delta_m.
    $$
    \end{itemize}
\end{itemize}
$\Box$\\

Now, we prove Theorem~~\ref{th-Ktildeixy}.
The points displayed in Lemma~\ref{lemma-discussion} entail that,
for $i\in E^{\rm o}$ and for $n\in\mathbb{N}^*$,
\begin{align}
\mathbb{P}_i\{X_1\in E^{\pm},\tilde{T}_n=0, X_n\in F\}=\;
& \mathbb{P}_i\{X_1\in E^-,\tau^+>n, X_n\in F\}
\nonumber\\
=\; & \sum_{j\in E^-} p_{ij}\,\mathbb{P}_j\{\tau^+>n-1,X_{n-1}\in F\},
\nonumber\\[-2ex]
\label{P-tildeTn}
\\[-2ex]
\mathbb{P}_i\{X_1\in E^{\pm},\tilde{T}_n=n,X_n\in F\}=\;
& \mathbb{P}_i\{X_1\in E^+,\tau^- >n,X_n\in F\}
\nonumber\\
=\; & \sum_{j\in E^+} p_{ij}\,\mathbb{P}_j\{\tau^- >n-1,X_{n-1}\in F\},
\nonumber
\end{align}
and for $i\in E^{\rm o}$ and $m\in\{1,\dots,n-1\}$,
\begin{equation}
\mathbb{P}_i\{X_1\in E^{\pm},\tilde{T}_n=m, X_n\in F\}
=\sum_{\ell=2}^n\sum_{j\in E^{\rm o}}
\mathbb{P}_i\{X_1\in E^{\pm},\tilde{\tau}^{\pm}=\ell,X_{\tilde{\tau}^{\pm}-1}=j,\tilde{T}_n=m,X_n\in F\}
\label{P-tildeTn-m}
\end{equation}
where, for $\ell\in\{2,\dots,n\}$,
\begin{align*}
\lqn{
\mathbb{P}_i\{X_1\in E^{\pm},\tilde{\tau}^{\pm}=\ell,X_{\tilde{\tau}^{\pm}-1}=j,\tilde{T}_n=m,X_n\in F\}}
\\[-3ex]
=\;& \mathbb{P}_i\{X_1\in E^+,\tilde{\tau}^{\pm}=\ell,X_{\tilde{\tau}^{\pm}-1}=j\}\,
\mathbb{P}_j\{X_1\in E^{\pm},\tilde{T}_{n-\ell+1}=m-\ell+1,X_{n-\ell+1}\in F\}
\\[1ex]
& +\mathbb{P}_i\{X_1\in E^-,\tilde{\tau}^{\pm}=\ell,X_{\tilde{\tau}^{\pm}-1}=j\}\,
\mathbb{P}_j\{X_1\in E^{\pm},\tilde{T}_{n-\ell+1}=m,X_{n-\ell+1}\in F\}.
\end{align*}
By~(\ref{P-tildeTn}) and by imitating the proof of~(\ref{eqK0xy}),
it can be easily seen that $\tilde{K}_i(x,0)$ and $\tilde{K}_i(0,y)$ can be written as~(\ref{tilde-K}).
Next, we compute $\tilde{K}_i(x,y)$. By~(\ref{P-tildeTn-m}),
\begin{align*}
\tilde{K}_i(x,y)=\;
& \tilde{K}_i(x,0)+\tilde{K}_i(0,y)
-\mathbb{P}_i\{X_0\in F,X_1\in E^{\pm}\}
+\sum_{m,n\in\mathbb{N}:\atop 1\le m\le n-1}
\mathbb{P}_i\{X_1\in E^{\pm},\tilde{T}_n=m, X_n\in F\}\,x^m y^{n-m}
\\
=\; & \tilde{K}_i(x,0)+\tilde{K}_i(0,y)
-\mathbbm{1}_F(i)\mathbb{P}_i\{X_1\in E^{\pm}\}
\\
& +\sum_{\ell=2}^{\infty}\sum_{j\in E^{\rm o}}
\mathbb{P}_i\{X_1\in E^+,\tilde{\tau}^{\pm}=\ell,X_{\tilde{\tau}^{\pm}-1}=j\}
\\
& \times \sum_{m,n\in\mathbb{N}:\atop\ell-1\le m\le n-1}
\mathbb{P}_j\{X_1\in E^{\pm},\tilde{T}_{n-\ell+1}=m-\ell+1, X_{n-\ell+1}\in F \}\,x^m y^{n-m}
\\
& +\sum_{\ell=2}^{\infty} \sum_{j\in E^{\rm o}}
\mathbb{P}_i\{X_1\in E^-,\tilde{\tau}^{\pm}=\ell,X_{\tilde{\tau}^{\pm}-1}=j\}
\\
& \times\sum_{m,n\in\mathbb{N}:\atop 1\le m\le n-\ell+1}
\mathbb{P}_j\{X_1\in E^{\pm},\tilde{T}_{n-\ell+1}=m, X_{n-\ell+1}\in F \}\,x^m y^{n-m}
\\
=\; & \tilde{K}_i(x,0)+\tilde{K}_i(0,y)-\mathbbm{1}_F(i)\mathbb{P}_i\{X_1\in E^{\pm}\}
\\
& +\sum_{\ell=2}^{\infty}\sum_{j\in E^{\rm o}} \mathbb{P}_i\{X_1\in E^+,
\tilde{\tau}^{\pm}=\ell,X_{\tilde{\tau}^{\pm}-1}=j\}\,x^{\ell-1}
\big(\tilde{K}_j(x,y)-\tilde{K}_j(x,0)\big)
\\
& +\sum_{\ell=2}^{\infty}
\sum_{j\in E^{\rm o}}\mathbb{P}_i\{X_1\in E^-,\tilde{\tau}^{\pm}=\ell,
X_{\tilde{\tau}^{\pm}-1}=j\}\,y^{\ell-1} \big(\tilde{K}_j(x,y)-\tilde{K}_j(0,y)\big)
\\
=\; & \tilde{K}_i(x,0)+\tilde{K}_i(0,y)
+\sum_{j\in E^{\rm o}} \big(\tilde{H}_{ij}^{\pm+}(x)
+\tilde{H}_{ij}^{\pm-}(y)\big)\tilde{K}_j(x,y)
\\
& -\sum_{j\in E^{\rm o}}\tilde{H}_{ij}^{\pm+}(x)\tilde{K}_j(x,0)
-\sum_{j\in E^{\rm o}}\tilde{H}_{ij}^{\pm-} (y)\tilde{K}_j(0,y)
-\mathbbm{1}_F(i)\mathbb{P}_i\{X_1\in E^{\pm}\},
\end{align*}
as claimed.
$\Box$

\subsection{Proof of Theorem \ref{th-gene-M=2}}\label{proof-th-gene-M=2}

\begin{itemize}
\item
\textsl{First step: writing out the matrices}

Suppose $F=\mathbb{Z}$ and set $\varpi=\pi_1+2\pi_2$,
$\mathbf{1}=\begin{pmatrix} 1 \\ 1\end{pmatrix}$,
$\mathbf{1\hspace{-.28em}l}=\begin{pmatrix} 1 \\ 1 \\ 1 \\ 1\end{pmatrix}$\!.
Below, we rewrite the matrices displayed in Example~\ref{example-M=2} which
can be simplified under the condition of symmetry $\pi_{-1}=\pi_1$, $\pi_{-2}=\pi_2$ and
$\Gamma_{-j}=\Gamma_j$ for any integer $j$:
$$
\mathbf{I}=\begin{pmatrix}1 & 0\\ 0 & 1\end{pmatrix}\!,\quad
\mathbf{I^{\boldsymbol{\pm}}}=\varpi\mathbf{I},\quad
\mathbf{1_F}=\mathbf{1},\quad
\mathbf{1_F^{\boldsymbol{\pm}}}=\varpi\mathbf{1},
$$
$$
\mathbf{P}=\begin{pmatrix} \pi_0 & \pi_1 \\ \pi_1 & \pi_0\end{pmatrix}\!,\quad
\mathbf{P_{\boldsymbol{-}}}=\begin{pmatrix} \pi_2 & \pi_1 \\ 0 & \pi_2\end{pmatrix}\!,\quad
\mathbf{P_{\boldsymbol{+}}}=\begin{pmatrix} \pi_2 & 0 \\ \pi_1 & \pi_2\end{pmatrix}\!,\quad
\mathbf{P_{\boldsymbol{\dag}}}=\begin{pmatrix} \pi_0 & \pi_1 & \pi_2 & 0 \\ \pi_1 & \pi_0 & \pi_1 & \pi_2\end{pmatrix}\!,\quad
$$
$$
\mathbf{\Gamma}(x)
=\begin{pmatrix} \Gamma_0(x) & \Gamma_1(x)\\ \Gamma_1(x) & \Gamma_0(x)\end{pmatrix}\!,\quad
\mathbf{\Gamma_{\boldsymbol{-}}}(x)
=\begin{pmatrix} \Gamma_2(x) & \Gamma_1(x)\\ \Gamma_3(x) & \Gamma_2(x)\end{pmatrix}\!,\quad
\mathbf{\Gamma_{\boldsymbol{+}}}(x)
=\begin{pmatrix} \Gamma_2(x) & \Gamma_3(x)\\ \Gamma_1(x) & \Gamma_2(x)\end{pmatrix}\!,
$$
$$
\mathbf{\Gamma_{\boldsymbol{=}}}(x)
=\begin{pmatrix} \Gamma_4(x) & \Gamma_3(x)\\ \Gamma_5(x) & \Gamma_4(x)\end{pmatrix}\!,\quad
\mathbf{\Gamma_{\boldsymbol{\ddag}}}(x)
=\begin{pmatrix} \Gamma_4(x) & \Gamma_5(x)\\ \Gamma_3(x) & \Gamma_4(x)\end{pmatrix}\!,
$$
$$
\mathbf{G}(x)=\mathbf{G_{\boldsymbol{+}}}(x)=\mathbf{G_{\boldsymbol{-}}}(x)
=\frac{1}{1-x}\,\mathbf{1},\quad
\mathbf{G_{\boldsymbol{\dag}}}(x)=\frac{1}{1-x}\,\mathbf{1\hspace{-.28em}l},
$$
$$
\mathbf{K}(x,y)=\begin{pmatrix} K_0(x,y) \\ K_1(x,y)\end{pmatrix}\!,\quad
\mathbf{\tilde{K}}(x,y)=\begin{pmatrix} \tilde{K}_0(x,y) \\ \tilde{K}_1(x,y)\end{pmatrix}\!.
$$
We observe that $\mathbf{P}$ and $\mathbf{\Gamma}(x)$
are symmetric and that $\mathbf{P_{\boldsymbol{+}}}$, $\mathbf{\Gamma_{\boldsymbol{+}}}(x)$,
$\mathbf{\Gamma_{\boldsymbol{\ddag}}}(x)$ are the transposes of $\mathbf{P_{\boldsymbol{-}}}$,
$\mathbf{\Gamma_{\boldsymbol{-}}}(x)$, $\mathbf{\Gamma_{\boldsymbol{=}}}(x)$, respectively.
Let us write the inverses of $\mathbf{\Gamma}(x)$ and
$\mathbf{I}-x\,\mathbf{P}$. For this, we introduce their
respective determinants:
$$
\delta(x)=(\pi_0^2-\pi_1^2)x^2-2\pi_0x+1,\quad
\Delta(x)=\Gamma_0(x)^2-\Gamma_1(x)^2.
$$
We have
$$
\mathbf{\Gamma}(x)^{-1}
=\frac{1}{\Delta(x)}\begin{pmatrix} \Gamma_0(x) & -\Gamma_1(x)\\
-\Gamma_1(x) & \Gamma_0(x)\end{pmatrix}\!,\quad
(\mathbf{I}-x\,\mathbf{P})^{-1}
=\frac{1}{\delta(x)}\begin{pmatrix} 1-\pi_0x & \pi_1x \\ \pi_1x & 1-\pi_0x\end{pmatrix}\!.
$$

\item
\textsl{Second step: deriving the generating matrix $\mathbf{K}$}

In view of~(\ref{matrixK-LR-complete}), we write $\mathbf{K}(x,y)$ as
$\mathbf{K}(x,y)=\mathbf{D}(x,y)^{-1}\mathbf{N}(x,y)$ where
\begin{align*}
\mathbf{D}(x,y)
&
=\mathbf{I} -x\left(\mathbf{P}+
\mathbf{D_{\boldsymbol{+}}}(x)+\mathbf{D_{\boldsymbol{-}}}(y)
\right)
\\
\mathbf{N}(x,y)
&
=\frac{1}{1-x}\big[\mathbf{I}-x\big(\mathbf{P}+
\mathbf{D_{\boldsymbol{+}}}(x)\big)\big]\mathbf{N_{\boldsymbol{-}}}(x)
+\frac{1}{1-y}\,\mathbf{N_{\boldsymbol{\dag}}}(y)-\mathbf{1}
\end{align*}
with
$$
\mathbf{D_{\boldsymbol{+}}}(x)
=\mathbf{P_{\boldsymbol{+}}}\mathbf{\Gamma_{\boldsymbol{-}}}(x)\mathbf{\Gamma}(x)^{-1},\quad
\mathbf{D_{\boldsymbol{-}}}(y)
=\mathbf{P_{\boldsymbol{-}}}\mathbf{\Gamma_{\boldsymbol{+}}}(y)\mathbf{\Gamma}(y)^{-1},
$$
$$
\mathbf{N_{\boldsymbol{-}}}(x)
=\mathbf{1}-\mathbf{\Gamma_{\boldsymbol{-}}}(x)
\mathbf{\Gamma}(x)^{-1}\mathbf{1},\quad
\mathbf{N_{\boldsymbol{\dag}}}(y)
=\mathbf{1}-y\,\mathbf{P_{\boldsymbol{\dag}}}\mathbf{1\hspace{-.28em}l}
-y\,\mathbf{D_{\boldsymbol{-}}}(y)\mathbf{1}.
$$
First, we calculate
\begin{align*}
\mathbf{D_{\boldsymbol{+}}}(x)
&
=\frac{1}{\Delta(x)}
\begin{pmatrix} \pi_2 &\! 0 \\ \pi_1 &\! \pi_2\end{pmatrix}\!\!
\begin{pmatrix} \Gamma_2(x) &\! \Gamma_1(x)\\ \Gamma_3(x) &\! \Gamma_2(x)\end{pmatrix}\!\!
\begin{pmatrix} \Gamma_0(x) &\!\!\! -\Gamma_1(x)\\ -\Gamma_1(x) &\!\!\! \Gamma_0(x)\end{pmatrix}
=\frac{1}{\Delta(x)}
\begin{pmatrix} D_{00}^+(x) & D_{01}^+(x) \\[1ex]
D_{10}^+(x) & D_{11}^+(x)
\end{pmatrix}
\vspace{-\baselineskip}
\end{align*}
with
\vspace{-.25\baselineskip}
\begin{align*}
D_{00}^+(x)
&
=\pi_2\big(\Gamma_0(x)\Gamma_2(x)-\Gamma_1(x)^2\big)
\\
D_{01}^+(x)
&
=\pi_2\big(\Gamma_0(x)\Gamma_1(x)-\Gamma_1(x)\Gamma_2(x)\big),
\\
D_{10}^+(x)
&
=\pi_1\big(\Gamma_0(x)\Gamma_2(x)-\Gamma_1(x)^2\big)
+\pi_2\big(\Gamma_0(x)\Gamma_3(x)-\Gamma_1(x)\Gamma_2(x)\big),
\\
D_{11}^+(x)
&
=\pi_1\big(\Gamma_0(x)\Gamma_1(x)-\Gamma_1(x)\Gamma_2(x)\big)
+\pi_2\big(\Gamma_0(x)\Gamma_2(x)-\Gamma_1(x)\Gamma_3(x)\big).
\end{align*}
Similarly,
\begin{align*}
\mathbf{D_{\boldsymbol{-}}}(y)
&
=\frac{1}{\Delta(y)}
\begin{pmatrix} \pi_2 &\! \pi_1 \\ 0 &\! \pi_2\end{pmatrix}\!\!
\begin{pmatrix} \Gamma_2(y) &\! \Gamma_3(y)\\ \Gamma_1(y) &\! \Gamma_2(y)\end{pmatrix}\!\!
\begin{pmatrix} \Gamma_0(y) &\!\!\! -\Gamma_1(y)\\ -\Gamma_1(y) &\!\!\! \Gamma_0(y)\end{pmatrix}
=\frac{1}{\Delta(y)}
\begin{pmatrix} D_{00}^-(y) & D_{01}^-(y) \\[1ex]
D_{10}^-(y) & D_{11}^-(y)
\end{pmatrix}
\vspace{-\baselineskip}
\end{align*}
with
\vspace{-.25\baselineskip}
\begin{align*}
D_{00}^-(y)
&
=\pi_1\big(\Gamma_0(y)\Gamma_1(y)-\Gamma_1(y)\Gamma_2(y)\big)
+\pi_2\big(\Gamma_0(y)\Gamma_2(y)-\Gamma_1(y)\Gamma_3(y)\big),
\\
D_{01}^-(y)
&
=\pi_1\big(\Gamma_0(y)\Gamma_2(y)-\Gamma_1(y)^2\big)
+\pi_2\big(\Gamma_0(y)\Gamma_3(y)-\Gamma_1(y)\Gamma_2(y)\big),
\\
D_{10}^-(y)
&
=\pi_2\big(\Gamma_0(y)\Gamma_1(y)-\Gamma_1(y)\Gamma_2(y)\big),
\\
D_{11}^-(y)
&
=\pi_2\big(\Gamma_0(y)\Gamma_2(y)-\Gamma_1(y)^2\big).
\end{align*}
With this at hand,
\begin{align}
\mathbf{D}(x,y)^{-1}
&=\frac{\Delta(x)\Delta(y)}{D(x,y)}\left(\begin{matrix}
(1-\pi_0x)\Delta(x)\Delta(y)-x\Delta(y)D_{11}^+(x)
-y\Delta(x)D_{11}^-(y)
\\
\pi_1x\Delta(x)\Delta(y)+x\Delta(y)D_{10}^+(x)
+y\Delta(x)D_{10}^-(y)
\end{matrix}\right.
\nonumber\\[1ex]
&\hspace{8em}\left.\begin{matrix}
\pi_1x\Delta(x)\Delta(y)+x\Delta(y)D_{01}^+(x)
+y\Delta(x)D_{01}^-(y)
\\
(1-\pi_0x)\Delta(x)\Delta(y)-x\Delta(y)D_{00}^+(x)
-y\Delta(x)D_{00}^-(y)
\end{matrix}\right)
\label{D}
\end{align}
where $D(x,y)$ is the determinant of the previous displayed matrix.

On the other hand, we have
\begin{align*}
\mathbf{N_{\boldsymbol{-}}}(x)
&
=\begin{pmatrix} 1 \\ 1\end{pmatrix}\!\!-\frac{1}{\Delta(x)}
\begin{pmatrix} \Gamma_2(x) & \Gamma_1(x)\\ \Gamma_3(x) & \Gamma_2(x)\end{pmatrix}\!\!
\begin{pmatrix} \Gamma_0(x) & -\Gamma_1(x)\\ -\Gamma_1(x) & \Gamma_0(x)\end{pmatrix}\!\!
\begin{pmatrix} 1 \\ 1\end{pmatrix}
=\frac{1}{\Delta(x)}
\begin{pmatrix} N_0^-(x) \\[1ex] N_1^-(x)
\end{pmatrix}
\vspace{-\baselineskip}
\end{align*}
with
\vspace{-.25\baselineskip}
\begin{align*}
N_0^-(x)
&
=\Delta(x)-\big(\Gamma_0(x)-\Gamma_1(x)\big)\big(\Gamma_1(x)+\Gamma_2(x)\big),
\\
N_1^-(x)
&
=\Delta(x)-\big(\Gamma_0(x)-\Gamma_1(x)\big)\big(\Gamma_2(x)+\Gamma_3(x)\big).
\end{align*}
Next, we calculate $\mathbf{N_{\boldsymbol{\dag}}}(y)$:
$$
\mathbf{N_{\boldsymbol{\dag}}}(y)=
\begin{pmatrix} 1 \\ 1\end{pmatrix}
-y\begin{pmatrix} \pi_0 & \pi_1 & \pi_2 & 0 \\ \pi_1 & \pi_0 & \pi_1 & \pi_2\end{pmatrix}
\!\!\begin{pmatrix} 1 \\ 1 \\ 1 \\ 1\end{pmatrix}
-\frac{y}{\Delta(y)}\begin{pmatrix} D_{00}^-(y) & D_{01}^-(y) \\[1ex]
D_{10}^-(y) & D_{11}^-(y)
\end{pmatrix}
\!\!\begin{pmatrix} 1 \\ 1\end{pmatrix}
=\frac{1}{\Delta(y)} \begin{pmatrix} N_0^{\dag}(y) \\[1ex] N_1^{\dag}(y)
\end{pmatrix}
\vspace{-\baselineskip}
$$
with
\vspace{-.25\baselineskip}
\begin{align*}
N_0^{\dag}(y)=\;
& \big(1-(1-\pi_1-\pi_2)y\big)\Delta(y)-yD_{00}^-(y)-yD_{01}^-(y)
\\
=\;
& \big(1-(1-\pi_1-\pi_2)y\big)\Delta(y)
-y\big(\Gamma_0(y)-\Gamma_1(y)\big)\big(\pi_1\Gamma_1(y)+(\pi_1+\pi_2)\Gamma_2(y)
+\pi_2\Gamma_3(y)\big),
\\
N_1^{\dag}(y)=\;
& \big(1-(1-\pi_2)y\big)\Delta(y)-yD_{10}^-(y)-yD_{11}^-(y)
\\
=\;
& \big(1-(1-\pi_2)y\big)\Delta(y)
-\pi_2y\big(\Gamma_0(y)-\Gamma_1(y)\big)\big(\Gamma_1(y)+\Gamma_2(y)\big).
\end{align*}

As a by-product, we obtain the following representation:
\begin{equation}\label{N}
\mathbf{N}(x,y)=\begin{pmatrix} N_0(x,y)\\ N_1(x,y) \end{pmatrix}
\vspace{-\baselineskip}
\end{equation}
with
\vspace{-.25\baselineskip}
\begin{align*}
N_0(x,y)=\;
& \frac{1}{(1-x)\Delta(x)^2}\left[\big((1-\pi_0x)\Delta(x)
-xD_{00}^+(x)\big) N_0^-(x)-x\big(\pi_1\Delta(x)+D_{01}^+(x)\big)N_1^-(x)\right]
\\
& +\frac{1}{(1-y)\Delta(y)}\,N_0^{\dag}(y)-1,
\\
N_1(x,y)=\;
& \frac{1}{(1-x)\Delta(x)^2}\left[-x\big(\pi_1\Delta(x)+D_{10}^+(x)\big)N_0^-(x)
+\big((1-\pi_0x)\Delta(x)-xD_{11}^+(x)\big)N_1^-(x)\right]
\\
& +\frac{1}{(1-y)\Delta(y)}\,N_1^{\dag}(y)-1.
\end{align*}

\item
\textsl{Third step: deriving the generating matrix $\mathbf{\tilde{K}}$}

In view of~(\ref{matrixK-tilde-LR-complete}), we write $\mathbf{\tilde{K}}(x,y)$ as
$\mathbf{\tilde{K}}(x,y)=\mathbf{\tilde{D}}(x,y)^{-1}\mathbf{\tilde{N}}(x,y)$
where
\begin{align*}
\mathbf{\tilde{D}}(x,y)=\;
& \mathbf{I}-\varpi x\,\mathbf{\tilde{D}_{\boldsymbol{+}}}(x)
-\varpi y\,\mathbf{\tilde{D}_{\boldsymbol{-}}}(y),
\\
\mathbf{\tilde{N}}(x,y)=\;
& \big[\mathbf{I}-\varpi x\,\mathbf{\tilde{D}_{\boldsymbol{+}}}(x)\big]
\big[\varpi\mathbf{1}+\frac{x}{1-x}\,\mathbf{\tilde{N}_{\boldsymbol{+}}}(x)\big]
+\big[\mathbf{I}-\varpi y\,\mathbf{\tilde{D}_{\boldsymbol{-}}}(y)\big]
\big[\varpi\mathbf{1}
+\frac{y}{1-y}\,\mathbf{\tilde{N}_{\boldsymbol{-}}}(y)\big]-\varpi\mathbf{1},
\vspace{-\baselineskip}
\end{align*}
with
\vspace{-.25\baselineskip}
$$
\mathbf{\tilde{D}_{\boldsymbol{+}}}(x)
=\mathbf{D_{\boldsymbol{+}}}(x)(\mathbf{I}-x\,\mathbf{P})^{-1},\quad
\mathbf{\tilde{D}_{\boldsymbol{-}}}(y)
=\mathbf{D_{\boldsymbol{-}}}(y)(\mathbf{I}-y\,\mathbf{P})^{-1},
$$
$$
\mathbf{\tilde{N}_{\boldsymbol{+}}}(x)
=\mathbf{P_{\boldsymbol{+}}} \big(\mathbf{1}
-\mathbf{\Gamma_{\boldsymbol{=}}}(x) \mathbf{\Gamma}(x)^{-1}\mathbf{1}\big),\quad
\mathbf{\tilde{N}_{\boldsymbol{-}}}(y)
=\mathbf{P_{\boldsymbol{-}}} \big(\mathbf{1}
-\mathbf{\Gamma_{\boldsymbol{\ddag}}}(y) \mathbf{\Gamma}(y)^{-1}\mathbf{1}\big).
$$
For instance,
$$
\mathbf{\tilde{D}_{\boldsymbol{+}}}(x)
=\frac{1}{\delta(x)\Delta(x)}\,\mathbf{D_{\boldsymbol{+}}}(x)
\begin{pmatrix} 1-\pi_0x &\!\!\! \pi_1x \\ \pi_1x &\!\!\! 1-\pi_0x\end{pmatrix}
=\frac{1}{\delta(x)\Delta(x)}\begin{pmatrix} \tilde{D}_{00}^+(x)
& \tilde{D}_{01}^+(x) \\[1ex]
\tilde{D}_{10}^+(x) & \tilde{D}_{11}^+(x)
\end{pmatrix}
\vspace{-\baselineskip}
$$
with
\vspace{-.25\baselineskip}
\begin{align*}
\tilde{D}_{00}^+(x)=\;
& \pi_2 \big[\pi_1x\,\Gamma_0(x)\Gamma_1(x)+(1-\pi_0x)\Gamma_0(x)\Gamma_2(x)
-(1-\pi_0x)\Gamma_1(x)^2-\pi_1x\,\Gamma_1(x)\Gamma_2(x)\big],
\\
\tilde{D}_{01}^+(x)=\;
& \pi_2\big[(1-\pi_0x)\Gamma_0(x)\Gamma_1(x)+\pi_2x\Gamma_0(x)\Gamma_2(x)
-\pi_1x\,\Gamma_1(x)^2-(1-\pi_0x)\Gamma_1(x)\Gamma_2(x)\big],
\\
\tilde{D}_{10}^+(x)=\;
& \pi_1^2x\,\Gamma_0(x)\Gamma_1(x)+\pi_1\big(1+(\pi_2-\pi_0)x\big)\Gamma_0(x)\Gamma_2(x)
+\pi_2(1-\pi_0x)\Gamma_0(x)\Gamma_3(x)
\\
& -\pi_1(1-\pi_0x)\Gamma_1(x)^2
-\big(\pi_2+(\pi_1^2-\pi_0\pi_2)x\big)\Gamma_1(x)\Gamma_2(x)-\pi_1\pi_2x\,\Gamma_1(x)\Gamma_3(x),
\\
\tilde{D}_{11}^+(x)=\;
& \pi_1(1-\pi_0x)\Gamma_0(x)\Gamma_1(x)
+\big(\pi_2+(\pi_1^2-\pi_0\pi_2)x\big)\Gamma_0(x)\Gamma_2(x)+\pi_1\pi_2x\,\Gamma_0(x)\Gamma_3(x)
\\
&-\pi_1^2x\,\Gamma_1(x)^2-\pi_1\big(1+(\pi_2-\pi_0)x\big)\Gamma_1(x)\Gamma_2(x)
-\pi_2(1-\pi_0x)\Gamma_1(x)\Gamma_3(x).
\end{align*}
Similarly, we obtain that
$$
\mathbf{\tilde{D}_{\boldsymbol{-}}}(y)
=\frac{1}{\delta(y)\Delta(y)}\begin{pmatrix} \tilde{D}_{00}^-(y)
& \tilde{D}_{01}^-(y) \\[1ex]
\tilde{D}_{10}^-(y) & \tilde{D}_{11}^-(y)
\end{pmatrix}
\vspace{-\baselineskip}
$$
with
\vspace{-.25\baselineskip}
\begin{align*}
\tilde{D}_{00}^-(y)=\;
& \pi_1(1-\pi_0y)\Gamma_0(y)\Gamma_1(y)+\big(\pi_2+(\pi_1^2-\pi_0\pi_2)y\big)\Gamma_0(y)\Gamma_2(y)
+\pi_1\pi_2y\,\Gamma_0(y)\Gamma_3(y)
\\
& -\pi_1^2y\,\Gamma_1(y)^2-\pi_1\big(1+(\pi_2-\pi_0)y\big)\Gamma_1(y)\Gamma_2(y)
-\pi_2(1-\pi_0y)\Gamma_1(y)\Gamma_3(y),
\\
\tilde{D}_{01}^-(y)=\;
& \pi_1^2y\,\Gamma_0(y)\Gamma_1(y)
+\pi_1\big(1+(\pi_2-\pi_0)y\big)\Gamma_0(y)\Gamma_2(y)+\pi_2(1-\pi_0y)\Gamma_0(y)\Gamma_3(y)
\\
& -\pi_1(1-\pi_0y)\Gamma_1(y)^2-\big(\pi_2+(\pi_1^2-\pi_0\pi_2)y\big)\Gamma_1(y)\Gamma_2(y)
-\pi_1\pi_2y\,\Gamma_1(y)\Gamma_3(y),
\\
\tilde{D}_{10}^-(y)=\;
& \pi_2\big[(1-\pi_0y)\Gamma_0(y)\Gamma_1(y)+\pi_1y\,\Gamma_0(y)\Gamma_2(y)
-\pi_1y\,\Gamma_1(y)^2-(1-\pi_0y)\Gamma_1(y)\Gamma_2(y)\big],
\\
\tilde{D}_{11}^-(y)=\;
& \pi_2\big[\pi_1y\,\Gamma_0(y)\Gamma_1(y)+(1-\pi_0y)\Gamma_0(y)\Gamma_2(y)
-(1-\pi_0y)\Gamma_1(y)^2-\pi_1y\,\Gamma_1(y)\Gamma_2(y)\big].
\end{align*}
With this at hand,
\begin{align}
\mathbf{\tilde{D}}(x,y)^{-1}=\;
& \frac{\delta(x)\delta(y)\Delta(x)\Delta(y)}{\tilde{D}(x,y)}\left(\begin{matrix}
\delta(x)\delta(y)\Delta(x)\Delta(y)-\varpi x\delta(y)\Delta(y)\tilde{D}_{00}^+(x)
-\varpi y\delta(x)\Delta(x)\tilde{D}_{00}^-(y)
\\
\varpi x\delta(y)\Delta(y)\tilde{D}_{01}^+(x)
+\varpi y\delta(x)\Delta(x)\tilde{D}_{01}^-(y)
\end{matrix}\right.
\nonumber\\[1ex]
&\hspace{5.5em}\left.\begin{matrix}
\varpi x\delta(y)\Delta(y)\tilde{D}_{10}^+(x)
+\varpi y\delta(x)\Delta(x)\tilde{D}_{10}^-(y)
\\
\delta(x)\delta(y)\Delta(x)\Delta(y)-\varpi x\delta(y)\Delta(y)\tilde{D}_{11}^+(x)
-\varpi y\delta(x)\Delta(x)\tilde{D}_{11}^-(y)
\end{matrix}\right)
\label{Dtilde}
\end{align}
where $\tilde{D}(x,y)$ is the determinant of the previous displayed matrix.

Concerning $\mathbf{\tilde{N}}(x)$,
\begin{align*}
\mathbf{\tilde{N}_{\boldsymbol{+}}}(x)=\;
& \begin{pmatrix} \pi_2 & 0 \\ \pi_1 & \pi_2\end{pmatrix}\!\!
\left[\begin{pmatrix} 1 \\ 1\end{pmatrix}
-\frac{1}{\Delta(x)}\begin{pmatrix} \Gamma_4(x) & \Gamma_3(x)\\
\Gamma_5(x) & \Gamma_4(x)\end{pmatrix}\!\!
\begin{pmatrix} \Gamma_0(x) & -\Gamma_1(x)\\ -\Gamma_1(x) & \Gamma_0(x)\end{pmatrix}\!\!
\begin{pmatrix} 1 \\ 1\end{pmatrix}\right]
=\frac{1}{\Delta(x)}\begin{pmatrix} \tilde{N}_0^+(x) \\[1ex] \tilde{N}_1^+(x) \end{pmatrix}
\vspace{-\baselineskip}
\end{align*}
with
\vspace{-.25\baselineskip}
\begin{align*}
\tilde{N}_0^+(x)=\;
& \pi_2\Delta(x)-\pi_2\big(\Gamma_0(x)-\Gamma_1(x)\big)\big(\Gamma_3(x)+\Gamma_4(x)\big),
\\
\tilde{N}_1^+(x)=\;
& (\pi_1+\pi_2)\Delta(x)-\big(\Gamma_0(x)-\Gamma_1(x)\big)
\big(\pi_1\Gamma_3(x)+(\pi_1+\pi_2)\Gamma_4(x)+\pi_2\Gamma_5(x)\big).
\end{align*}
Similarly,
$$
\mathbf{\tilde{N}_{\boldsymbol{-}}}(y)=\frac{1}{\Delta(y)}
\begin{pmatrix} \tilde{N}_0^-(y) \\[1ex] \tilde{N}_1^-(y) \end{pmatrix}
\vspace{-\baselineskip}
$$
with
\vspace{-.25\baselineskip}
\begin{align*}
\tilde{N}_0^-(y)=\;
& (\pi_1+\pi_2)\Delta(y)-\big(\Gamma_0(y)-\Gamma_1(y)\big)
\big(\pi_1\Gamma_3(y)+(\pi_1+\pi_2)\Gamma_4(y)+\pi_2\Gamma_5(y)\big),
\\
\tilde{N}_1^-(y)=\;
& \pi_2\Delta(y)-\pi_2\big(\Gamma_0(y)-\Gamma_1(y)\big)\big(\Gamma_3(y)+\Gamma_4(y)\big).
\end{align*}
As a by-product, we obtain the following representation:
\begin{equation}
\mathbf{\tilde{N}}(x,y)=\begin{pmatrix} \tilde{N}_0(x,y)\\ \tilde{N}_1(x,y) \end{pmatrix}
\label{Ntilde}
\vspace{-\baselineskip}
\end{equation}
with
\vspace{-.25\baselineskip}
\begin{align*}
\tilde{N}_0(x,y)=\;
& \Big[1-\frac{\varpi x}{\delta(x)\Delta(x)}\,\tilde{D}_{00}^+(x)\Big]\!
\Big[\varpi+\frac{x}{(1-x)\Delta(x)}\,\tilde{N}_0^+(x)\Big]
\\
& -\frac{\varpi x}{\delta(x)\Delta(x)}\,\tilde{D}_{01}^+(x)
\Big[\varpi+\frac{x}{(1-x)\Delta(x)}\,\tilde{N}_1^+(x)\Big]
\\
& +\Big[1-\frac{\varpi y}{\delta(y)\Delta(y)}\,\tilde{D}_{00}^-(y)\Big]\!
\Big[\varpi+\frac{y}{(1-y)\Delta(y)}\,\tilde{N}_0^-(y)\Big]
\\
& -\frac{\varpi y}{\delta(y)\Delta(y)}\,\tilde{D}_{01}^-(y)
\Big[\varpi+\frac{y}{(1-y)\Delta(y)}\,\tilde{N}_1^-(y)\Big]
-\varpi,
\\
\tilde{N}_1(x,y)=\;
& -\frac{\varpi x}{\delta(x)\Delta(x)}\,\tilde{D}_{10}^+(x)
\Big[\varpi+\frac{x}{(1-x)\Delta(x)}\,\tilde{N}_0^+(x)\Big]
\\
& +\Big[1-\frac{\varpi x}{\delta(x)\Delta(x)}\,\tilde{D}_{11}^+(x)\Big]\!
\Big[\varpi+\frac{x}{(1-x)\Delta(x)}\,\tilde{N}_1^+(x)\Big]
\\
& -\frac{\varpi y}{\delta(y)\Delta(y)}\,\tilde{D}_{10}^-(y)
\Big[\varpi+\frac{y}{(1-y)\Delta(y)}\,\tilde{N}_0^-(y)\Big]
\\
& +\Big[1-\frac{\varpi y}{\delta(y)\Delta(y)}\,\tilde{D}_{11}^-(y)\Big]\!
\Big[\varpi+\frac{y}{(1-y)\Delta(y)}\,\tilde{N}_1^-(y)\Big]-\varpi.
\end{align*}
\end{itemize}
$\Box$

\vspace{1cm}

\noindent \textsc{Acknowledgements.} The authors thank the anonymous referees
who provide many constructive suggestions for improving the presentation
of the paper.
\\


\end{document}